\documentclass[11pt,a4paper]{article}
\usepackage[utf8]{inputenc}
\usepackage{amsmath}
\usepackage{amsfonts}
\usepackage{amssymb}
\usepackage{makeidx}
\usepackage{graphicx}
\usepackage{mathabx} %benötigt für schreibweise der rankine-hugoniot bedingung
\usepackage{dsfont} %benötigt für charakteristische funktion (1-funktion)
\usepackage[left=2cm,right=2cm,top=2cm,bottom=2cm]{geometry}
\usepackage{color}
\usepackage{verbatim} %füge kommentare ein durch \begin{comment} \end{comment}
\usepackage{tocbibind} %packt literaturverzeichnis in table of contents (allerdings auch table of contents selbst)
\usepackage{hyperref} %um hyperlinks einzufügen
\usepackage{mathtools} %für inklusionspfeil

\newtheorem{satz}{Lemma}[section]
\newtheorem{proposition}{Proposition}[section]

\newtheorem{remark}{Remark}[section]
\newtheorem{definition}{Definition}[section]
\newtheorem{theorem}{Theorem}[section]
\begin{document}

\setlength{\parindent}{0pt} %damit nicht jeder neue absatz eingerückt ist eingerückt ist

\begin{center}
\Large\textbf{DISSIPATIVE SOLUTIONS TO THE MODEL OF A GENERAL COMPRESSIBLE VISCOUS FLUID WITH THE COULOMB FRICTION LAW BOUNDARY CONDITION} \\[10mm]
\large{\v{S}ÁRKA NE\v{C}ASOVÁ}$^1$, \large{JUSTYNA OGORZAŁY}$^{1,3}$, \large{JAN SCHERZ}$^{1,2}$ \\[10mm]
\end{center}

\begin{itemize}
\item[$^1$] Institute of Mathematics of the Czech Academy of Sciences, \v{Z}itná 25, Prague 1, 11567
\item[$^2$] Department of Mathematical Analysis, Faculty of Mathematics and Physics, Charles University in Prague, Sokolovská 83, Prague 8, 18675
\item[$^3$] Institute of Mathematics and Cryptology, Military University of Technology, ul. Gen. Sylwestra Kaliskiego 2,
00-908 Warsaw 
\end{itemize}

\begin{center}
\Large\textbf{Abstract} \\[4mm]
\end{center}
      
We study a model of a general compressible viscous fluid subject to the Coulomb friction law boundary condition. For this model, we introduce a dissipative formulation and prove the existence of dissipative solutions. The proof of this result consists of a three level approximation method: A Galerkin approximation, the classical parabolic regularization of the continuity equation as well as convex regularizations of the potential generating the viscouss stress and the boundary terms incorporating the Coulomb friction law into the dissipative formulation. This approach combines the techniques already known from the proof of the existence of dissipative solutions to a model of general compressible viscous fluids under inflow-outflow boundary conditions as well as the proof of the existence of weak solution to the incompressible Navier-Stokes equations under the Coulomb friction law boundary condition. It is the first time that this type of boundary condition is considered in the case of compressible flow.

\section{Introduction}

In this article we prove the existence of dissipative solutions to a model of a general compressible viscous fluid subject to the Coulomb friction law boundary condition.

The theory of the existence of weak solutions to the Navier-Stokes equations is well developed. For an introduction to the topic we refer to the books \cite{lions} - in the case of an incompressible fluid - and \cite{lions2} - in the case of a compressible fluid. The existence of weak solutions is proved in these books in the cases of the no-slip boundary condition - meaning that the fluid velocity vanishes on the boundary of the considered domain -, periodic boundary conditions as well as in the case of the fluid occupying the whole space. In the setting of the compressible Navier-Stokes equations, the mathematical analysis is commonly focused on the isentropic case, meaning that the pressure is assumed to be the $\gamma$-th power of the density for a so-called adiabatic constant $\gamma \geq 1$. While in \cite{lions2} this constant is subject to the restriction $\gamma \geq \frac{9}{5}$ (in a three dimensional spatial domain), \,{ an extension to a more physical situation - the case $\gamma > \frac{3}{2}$,  which is the current state of the art - and the no-slip boundary condition is given in \cite{fnp} (for more details, see \cite{ns}). Moreover, we  refer to \cite{fn} for the proof of the existence of weak solutions to the full system of a compressible fluid with heat conductivity.} This result treats, in addition to the no-slip boundary condition, also the complete slip boundary condition, under which the normal component of the velocity vanishes on the domain boundary.

In the recent years, also rather uncommon boundary or some time "nonstandard" conditions have started to move towards the center of attention in the study of fluid dynamics. This is motivated by the hope of finding suitable ways for the description of fluid motion which do not lead to certain paradoxical phenoma occurring under the more classical boundary conditions; an example being the impossibility of a contact between rigid bodies moving freely through a linearly viscous fluid, cf.\ \cite{HES,HIL}. See also work of Bucur et al, in  the case of rough boundary \cite{BFN} or nonstandard boundary condition by Bellout et al. \cite{BNP}.

One specific example for a boundary condition which leads to, in at least some aspects, more physically reasonable results is the Navier-slip boundary condition, see e.g.\ \cite{MOF}. A proof of the existence of weak solutions to the model of the motion of a rigid body immersed in a compressible fluid under Navier-slip boundary and interface conditions can be found in \cite{compressiblecase}. Two further boundary conditions, which model the motion of fluids slipping on the boundary if the tangential component of the stress tensor is sufficiently large, are the slip boundary condition of friction type and the Coulomb friction law boundary condition. The first one of these conditions was introduced in \cite{fu,fuj2} for the setting of the stationary Stokes and Navier-Stokes equations. For the evolutionary incompressible Navier-Stokes equations, the existence of (strong) solutions under this boundary condition was proved in \cite{kas2} globally in time in two dimensional spatial domains and locally in time in three dimensional domains. Additionally, in \cite{slipbcoffrictiontype}, the global in time existence of weak solutions to the compressible Navier-Stokes equations was achieved under this condition. The existence of weak solutions to the incompressible Navier-Stokes equations subject to the Coulomb friction law boundary condition was proved in \cite{bmt}. Moreover, for the same problem with an additional rigid body moving through the fluid, the local in time existence of weak solutions was obtained in \cite{bmt2}.

Results on the existence of solutions to models of compressible fluids satisfying the Coulomb friction law boundary condition are still missing. This is where the present article comes in. However, instead of the classical compressible Navier-Stokes equations, we are interested in a more general model. More specifically, we consider the model which was analyzed in \cite{afn} in the setting of inflow-outflow boundary conditions. In this model, the viscous stress is generated by a suitable potential via a general implicit rheological law (cf.\ the relation \eqref{sdef} below) instead of given by the classical linear Newton's rheological law. Moreover, while the pressure is chosen such that its potential is still bounded from below by the $\gamma$-th power of the density, it is not restricted to the isentropic case and arbitrarily small values $\gamma > 1$ are taken into account.

As the notion of weak solutions is not eligible under these loosened restrictions, the article \cite{afn} instead studies the existence (as well as an associated weak-strong uniqueness property) of dissipative solutions to the considered system. In the setting of dissipative solutions, the Cauchy stress is supplemented by a so-called turbulent component, causing an error term in the momentum equation. In the case that this error term vanishes, the dissipative formulation turns into the standard weak formulation of the problem. Before the article \cite{afn}, the notion of dissipative solutions was already studied for the compressible Euler equations in \cite{bfh} and for general models of incompressible fluids in \cite{af}. We also mention the article \cite{b}, wherein the existence result for dissipative solutions in the compressible setting under inflow-outflow boundary conditions in \cite{afn} was generalized to the case of a linear pressure, i.e.\ the isentropic case with the adiabatic constant $\gamma = 1$.

In the present article we introduce a notion of dissipative solutions to the system from \cite{afn} under the Coulomb friction law boundary condition and we prove their existence. Our proof follows by and large the three level approximation method used for the proof of the existence of dissipative solutions in the case of inflow-outflow boundary conditions given in \cite{afn}. The three levels in this procedure consist of a Galerkin approximation, the classical parabolic regularization of the continuity equation (see \cite[Section 7.6]{ns}) as well as a suitable convex regularization of the potential $F$ determining the viscous stress. The greatest novelty in our proof is the combination of this approach with the approach used in \cite{bmt} for the proof of the existence of weak solutions to the incompressible Navier-Stokes equations under the Coulomb friction law boundary condition. More specifically, in parallel to the regularization of the potential $F$, we construct a convex regularization $j_\delta (u)$ of the absolute value of the velocity field $u$ and supplement the momentum equation by an additional boundary integral containing the gradient of $j_\delta$, see \eqref{m22}. Under exploitation of the convexity of $j_\delta$, the approximate momentum equation can then be turned into the desired inequality containing the boundary integrals expressing the Coulomb friction law boundary condition as in our dissipative formulation, cf.\ \eqref{momineq}.

A further difficulty arises in the limit passage in the regularization of the continuity equation: Here, in order to correctly identify the limits of the boundary integrals in the momentum equation, we require strong convergence of the velocity field. This, in turn, can be achieved under the condition of strong convergence of the density. Classically, strong convergence of the density is achieved by proving a so-called effective viscous flux identity, which, however, is not possible in our case since we still find ourselves in a Galerkin approximation at this point. Nonetheless, it turns out that the Galerkin regularity of the velocity field suffices to prove strong convergence of the density via a comparison between the renormalized continuity equation before and after the limit passage, even without an effective viscous flux identity at hand.

Finally, in the limit passage in the Galerkin method, we are less lucky: It does not appear to be possible to deduce strong convergence of the velocity field at this stage anymore. {\it This leads to one of the boundary integrals in the momentum equation remaining unidentified in our final system. As a compensation, however, we are able to at least identify this term in the sense of Young measures by transferring the classical theory of Young measures to the setting of boundary integrals,} cf.\ Lemma \ref{youngmeasures} in the Appendix.

The article is outlined as follows: In Section \ref{model} we present the considered model in its classical formulation. In Section \ref{dissform} we introduce our associated disspative formulation as well as our main result. In Section \ref{sub1}, we introduce the approximation method used for the proof of this result. Subsequently, in Section \ref{approxexist}, we solve the approximate system and finally, in Sections \ref{reglim}--\ref{gallim}, we conclude the proof by passing to the limit in all levels of the approximation.

\section{Model} \label{model}

We study a model of the motion of a general compressible viscous fluid subject to the Coulomb friction law boundary condition. Let $T>0$. The fluid is assumed to be contained in a bounded domain $\Omega \subset \mathbb{R}^3$ with boundary $\Gamma := \partial \Omega$ and outward unit normal vector $\operatorname{n}$. The motion of the fluid is described by its density $\rho: (0,T)\times \Omega \rightarrow \mathbb{R}$ and its velocity $u:(0,T)\times \Omega \rightarrow \mathbb{R}^3$, the evolution of which is determined by the system of partial differential equations
\begin{align}
\partial_t \rho + \nabla \cdot \left(\rho u \right) = 0 \quad \quad \quad \quad \quad \quad \quad \quad \quad &\text{in } (0,T)\times \Omega, \label{classiccont} \\
\partial_t \left( \rho u \right) + \nabla \cdot \left( \rho u \otimes u \right) = \nabla \cdot \mathbb{T} + \rho f \quad \quad \quad \quad \quad &\text{in } (0,T)\times \Omega, \label{classicmom} \\
\rho(0)= \rho_0,\quad (\rho u)(0) = (\rho u)_0 \quad \quad \quad \quad \quad \ \quad \quad &\text{in } \Omega, \label{initialconds} \\
u \cdot \operatorname{n} = 0 \quad \quad \quad \quad \quad \quad \quad \quad \quad &\text{on } (0,T)\times \Gamma, \label{boundcond1} \\
\left( \mathbb{T}\operatorname{n} \right)_\tau \cdot y \geq g \left| u_\tau \right| - g \left| u_\tau + y \right| \quad \quad &\text{on } (0,T)\times \Gamma \quad \forall y \in \mathbb{R}^3. \label{coulombbc}
\end{align}

In this system, the relations \eqref{classiccont} and \eqref{classicmom} constitute the continuity equation and the momentum equation. While the continuity equation does not differ from the corresponding relation in the classical compressible Navier-Stokes equations, the momentum equation is here allowed to take a more general form: Namely, the Cauchy stress $\mathbb{T}$ is assumed to be given as $\mathbb{T} = \mathbb{S} - p \operatorname{Id}$, where the viscous stress $\mathbb{S}$ is generated by a suitable convex potential $F:\mathbb{R}^{3\times 3}_{\operatorname{sym}} \rightarrow [0,\infty)$ via the implicit rheological law
\begin{align}
\mathbb{D}u : \mathbb{S} = F\left(\mathbb{D}u \right) + F^*(\mathbb{S}),\quad \quad \mathbb{D}u := \frac{1}{2} \nabla u + \frac{1}{2} \left( \nabla u \right)^T \label{sdef}
\end{align}

for the convex conjugate $F^*:\mathbb{R}^{3\times 3}_{\operatorname{sym}} \rightarrow [0,\infty]$ of $F$ and the pressure $p= p(\rho):[0,\infty) \rightarrow [0,\infty)$ is a function of the density $\rho$. In the classical compressible Navier-Stokes equations, the potential $F$ is given in the form
\begin{align}
F\left( \mathbb{D}u \right) = \frac{\mu}{2} \left| \mathbb{D}u \right|^2 + \frac{\lambda}{2} \left| \nabla \cdot u \right|^2 \label{classpot}
\end{align}

for some constants $\mu > 0$ and $\lambda \geq - \frac{2}{3}\mu$, which corresponds to the viscous stress given by the classical Newton's rheological law
\begin{align}
\mathbb{S} = \mu \mathbb{D}u + \lambda \left( \nabla \cdot u \right) \operatorname{Id}. \nonumber
\end{align}

The pressure $p$ in the mathematical analysis of the classical compressible Navier-Stokes equations is typically assumed to obey the isentropic constitutive relation
\begin{align}
p(\rho) = a\rho^\gamma \label{classpress}
\end{align}

for some constants $a>0$ and $\gamma > \frac{3}{2}$. In the present article, while we do impose certain restrictions on $F$ and $p$, these functions are not restricted to the cases \eqref{classpot} and \eqref{classpress}. The specific conditions we assume for $F$ and $p$ are the same ones as in the analysis of general compressible fluids in \cite{afn} and are detailed in Section \ref{dissform} below. The quantity $f:(0,T)\times \Omega \rightarrow \mathbb{R}^3$ in the momentum equation \eqref{classicmom} constitutes a given external forcing term. The relations \eqref{initialconds} in our model constitute the initial conditions with given initial data $\rho_0:\Omega \rightarrow \mathbb{R}$ and $(\rho u)_0:\Omega \rightarrow \mathbb{R}^3$. Finally, the relations \eqref{boundcond1} and \eqref{coulombbc} represent the Coulomb friction law boundary condition. Therein, the notation $v_\tau := v - (v\cdot \operatorname{n}) \operatorname{n}$ represents the tangential component of arbitrary vectors $v \in \mathbb{R}^3$ and the given function $g: (0,T)\times \Gamma \rightarrow [0,\infty)$ denotes the threshold of slippage.

\section{Dissipative solutions and main result} \label{dissform}

In this section we introduce our definition of dissipative solutions to the problem \eqref{classiccont}--\eqref{coulombbc} and state the main result. Moreover, we introduce some additional notation and structural restrictions on the potential $F$ and the pressure $p$. In addition to the standard notation for the Lebesgue and Sobolev spaces we use the notation
\begin{align}
W_{\operatorname{n}}^{k,p}\left(\Omega;\mathbb{R}^3\right) := \left\lbrace \phi \in W^{k,p}\left(\Omega;\mathbb{R}^3\right): \left. \phi \cdot \operatorname{n}\right|_\Gamma = 0 \right\rbrace \quad \quad \text{for } k=1,2,\ 1 \leq p \leq \infty \nonumber
\end{align}

for the spaces of Sobolev functions the normal component of which vanishes on the boundary $\Gamma$ of $\partial \Omega$. By $\mathcal{M}(\overline{\Omega})$ we denote the space of all Borel measures on $\overline{\Omega}$ and by $\mathcal{M}^+(\overline{\Omega})$ its subspace of all positive Borel measures on $\overline{\Omega}$. Further we generalize the latter space to a space of tensor-valued Borel measures,
\begin{align}
\mathcal{M}^+\left(\overline{\Omega}; \mathbb{R}^{3\times 3}_{\operatorname{sym}} \right) := \left\lbrace \mathcal{R} = \left( \mathcal{R}_{ij} \right)_{i,j=1}^3:\ \mathcal{R}_{ij} = \mathcal{R}_{ji} \in \mathcal{M}(\overline{\Omega}), \ \mathcal{R}: (\xi \otimes \xi) \in \mathcal{M}^+\left(\overline{\Omega}\right) \quad \forall \xi \in \mathbb{R}^3 \right\rbrace. \nonumber
\end{align}

For the potential $F$ we assume that
\begin{align}
F: \mathbb{R}^{3\times 3}_{\text{sym}} \rightarrow [0,\infty) \quad \text{is proper convex and lower semicontinuous with } F(0)=0. \label{fcond1}
\end{align} 

Additionally, for any $R>0$, we suppose the existence of a convex and increasing function $A_R:[0,\infty)\rightarrow [0,\infty)$ satisfying $A_R(0)=0$ as well as the $\Delta_2^2$-condition
\begin{align}
\alpha_1 A_R(z) \leq A_R(2z) \leq \alpha_2A_R(z) \quad \quad \text{for all } z \in [0,\infty) \nonumber
\end{align}

and two constants $\alpha_1 > 2$, $\alpha_2 < \infty$, such that
\begin{align}
F\left(D + Q \right) -F(D) - S:Q \geq A_r \left( \left| Q - \frac{1}{3} \operatorname{tr} \left[ Q \right] \operatorname{Id} \right| \right) \label{fcond2}
\end{align}

for all $D,Q,S \in \mathbb{R}^{3\times 3}_{\text{sym}}$ with $|D|\leq R$ and  $S \in \partial F(Q)$, where $\partial F(Q)$ denotes the subdifferential of $F$ in $Q$. The assumptions \eqref{fcond1} and \eqref{fcond2} imply several further useful properties for the potential $F$:

\begin{remark} \label{fprops}
The convex potential $F: \mathbb{R}^{3\times 3}_{\operatorname{sym}} \rightarrow [0,\infty)$ is in particular continuous, see for example \cite[Corollary 10.1.1]{rockafellar}. Furthermore, due to the assumptions \eqref{fcond1}, the convex conjugate $F^*:\mathbb{R}^{3\times 3}_{\operatorname{sym}} \rightarrow [0,\infty]$ of $F$ is also convex and lower semicontinuous as well as superlinear, cf.\ \cite[Proposition 2.1]{b}. The convexity and the lower semicontinuity of $F$ and $F^*$ further imply that both these functions are weakly lower semicontinuous. Moreover, in view of the proper convexity and the lower semicontinuity of $F$, the implicit rheological law \eqref{sdef} can be stated equivalently as the requirement of the tensor $\mathbb{S}$ being a subgradient of $F$ in $\mathbb{D}u$,
\begin{align}
\mathbb{S} \in \partial F \left( \mathbb{D}u \right) \quad \Leftrightarrow \quad \mathbb{D}u \in \partial F^* \left(\mathbb{S} \right) \quad \Leftrightarrow \quad \mathbb{D}u : \mathbb{S} = F\left(\mathbb{D}u \right) + F^*(\mathbb{S}), \label{fenchineq}
\end{align}
cf.\ \cite[Theorem 3.32]{rindler}. Finally, the estimate \eqref{fcond2} allows us to find constants $\mu >0$, $q>1$, such that
\begin{align}
F(D) \geqslant \mu \Big \vert D - \frac{1}{3} \operatorname{tr} \left[ D\right] \operatorname{Id} \Big \vert^q \quad \forall D \in \mathbb{R}^{3\times 3}_{\operatorname{sym}} \quad \text{with} \quad \vert D \vert >1, \label{fprop}
\end{align}

cf.\ \cite[Section 2.1.2]{afn}.
\end{remark}

For the pressure $p$ we suppose that
\begin{align}
p \in C\left([0,\infty) \right) \bigcap C^2\left(0,\infty \right),\quad \quad p(0)=0,\quad \quad p'(\rho)>0 \quad \forall \rho > 0. \label{pcond1}
\end{align}

Furthermore, we denote by $P$ the pressure potential associated to $p$, defined by
\begin{align}
P(\rho) = \rho \int_1^\rho \frac{p(z)}{z^2}\ dz \label{pcond2}
\end{align}

and we assume the existence of two constants $\underline{a}, \overline{a}>0$ such that the functions
\begin{align}
P(\cdot) - \underline{a}p(\cdot) \quad \text{and} \quad \overline{a}p(\cdot) - P(\cdot) \quad \text{are convex.} \label{pcond3}
\end{align}

Again, the assumptions \eqref{pcond1}--\eqref{pcond3} have several useful consequences for $p$ and $P$:
\begin{remark} \label{pprops}
The pressure potential $P$, which satisfies the relation
\begin{align}
P'(\rho)\rho - P(\rho) = p(\rho), \label{Pid}
\end{align}

obeys the coercivity condition
\begin{align}
P(\rho) \geq a\rho^\gamma \quad \forall \rho \geq 1 \label{coercP}
\end{align}

for some constants $a>0$ and $\gamma>1$, cf.\ \cite[Section 2.1.1]{afn}. Moreover, the pressure $p$ satisfies the estimate
\begin{align}
p''(\rho) \geq \frac{P''(\rho)}{\overline{a}} = \frac{p'(\rho)}{\overline{a}\rho} > 0 \quad \forall \rho > 0 \label{convexity}
\end{align}

and is thus in particular convex, cf.\ \cite[Section 2.1.1]{afn}. Finally, due to its $C^2(0,\infty)$-regularity, the pressure is locally Lipschitz-continuous in $(0,\infty)$.
\end{remark}

We now present our full definition of dissipative solutions to the problem \eqref{classiccont}--\eqref{coulombbc}, followed up by a brief explanation afterwards.
\begin{definition} \label{disssol}
Let $T > 0$ and let $\Omega \subset \mathbb{R}^3$ be a bounded domain with boundary $\Gamma = \partial \Omega$. Let the data
\begin{align}
f \in L^\infty((0,T)\times \Omega;\mathbb{R}^3),\quad \quad 0 \leq g \in L^\infty((0,T)\times \Gamma;\mathbb{R}),\quad \quad F:\mathbb{R}^{d \times d}_{\operatorname{sym}} \to [0, \infty),\quad \quad p:[0,\infty) \rightarrow \mathbb{R} \nonumber
\end{align}

satisfy the conditions \eqref{fcond1}, \eqref{fcond2} and \eqref{pcond1}--\eqref{pcond3}. Moreover, consider initial data
\begin{align}
\rho_0 \in L^\gamma(\Omega;\mathbb{R}),\quad \rho_0 \geqslant 0,\quad (\rho u)_0 \in L^{\frac{2 \gamma}{\gamma + 1}}\left(\Omega; \mathbb{R}^3\right),\quad \int_{\Omega}^{}\, \left( \frac{1}{2}\, \frac{\vert (\rho u)_0 \vert^2}{\rho_0} + P(\rho_0) \right)\, dx < \infty, \nonumber
\end{align}

where $P$ denotes the pressure potential defined in \eqref{pcond2} and $\gamma>1$ is chosen according to the estimate \eqref{coercP}. Then the problem \eqref{classiccont}--\eqref{coulombbc} is said to admit a dissipative solution 
\begin{align}
0 \leq \rho \in C_{\operatorname{weak}} \left([0,T];L^\gamma (\Omega;\mathbb{R}) \right),\quad \quad u \in L^q \left(0,T;W_{\operatorname{n}}^{1,q}(\Omega;\mathbb{R}^3) \right), \label{regularity1}
\end{align}
where $q>1$ is determined by the estimate \eqref{fprop}, if there exist
%\begin{align}
%\rho u \in C_{\operatorname{weak}} \left([0,T];L^\frac{2\gamma}{\gamma + 1}\left(\Omega; \mathbb{R}^3 \right) \right), \label{regularity2}
%\end{align}
%if there exist
\begin{align}
\mathbb{S} \in L^1\left((0,T)\times \Omega; \mathbb{R}_{\text{sym}}^{3\times 3} \right), \quad \mathcal{E} \in L^\infty \left((0,T; \mathcal{M}^+ \left(\overline{\Omega}; \mathbb{R}\right)\right),\quad \mathcal{R} \in L^\infty\left(0,T;\mathcal{M}^+ \left(\overline{\Omega}; \mathbb{R}^{3 \times 3}_{\text{sym}}\right)\right), \label{additionalterms}
\end{align}

where $\mathcal{E}$ and $\mathcal{R}$ are related via the compatibility condition
\begin{align}
\underline{d}\mathcal{E} \leq \operatorname{tr} \left[\mathcal{R} \right] \leq \overline{d}\mathcal{E} \label{compcon}
\end{align}

for some constants $0 < \underline{d} \leq \overline{d}< \infty$, and if there exists a family $\{\nu_{(t,y)} \}_{(t,y) \in (0,T)\times \Gamma}$ of probability measures $\nu_{(t,y)}$ on $\mathbb{R}^3$ depending measurably on $(t,y)$ as well as, for any $\phi \in \mathcal{D}([0,T]\times \overline{\Omega})$, a function
\begin{align}
\overline{|\phi + u|} \in L^q \left((0,T)\times \Gamma \right),\quad \quad \overline{|\phi + u|}(t,y) = \int_{\mathbb{R}^3} |\phi\left(t,y\right) + z|\ d\nu_{(t,y)}(z) \quad \text{for a.a. } (t,y) \in (0,T) \times \Gamma \label{boundaryid}
\end{align}

depending on $\phi$ continuously with respect to the $L^q((0,T)\times \Gamma)$-norm, such that the continuity equation
\begin{align}
\partial_t \rho + \nabla \cdot \left( \rho u \right) = 0 \quad \text{in } \mathcal{D}'\left((0,T)\times \Omega \right) \label{conteq}
\end{align}

holds true, the momentum inequality
\begin{align}
-\int_0^T \int_\Omega \rho u \cdot \partial_t \phi \ dxdt - \int_\Omega \left(\rho u\right)_0\phi(0)\ dx \geq& \int_0^T \int _\Omega \left( \rho u \otimes u \right): \nabla \phi- \mathbb{S}: \nabla \phi + p \left(\rho \right) \nabla \cdot \phi + \rho f\cdot \phi \ dxdt \nonumber \\
&+ \int_0^T \int_\Omega \nabla \phi :\ d\mathcal{R}(t)dt + \int_0^T \int_{\Gamma} g \left|u \right| - g \overline{\left| \phi + u \right|} \ d\Gamma dt \label{momineq}
\end{align}

is satisfied for all $\phi \in \mathcal{D}([0,T)\times \overline{\Omega})$ with $\phi \cdot \operatorname{n}|_{\Gamma} = 0$ and the energy inequality
\begin{align}
&\int_\Omega \frac{1}{2} \rho (\tau)|u(\tau)|^2 + P \left(\rho (\tau) \right)\ dx + \int_{\overline{\Omega}}1\ d\mathcal{E}(\tau) + \int_0^\tau \int_\Omega F\left(\mathbb{D} u\right) + F^* \left(\mathbb{S}\right)\ dxdt + \int_0^\tau \int_{\Gamma} g \left|u\right|\ d\Gamma dt \nonumber \\
\leq&\int_\Omega \frac{1}{2} \rho (\tau)|u(\tau)|^2 + P \left(\rho (\tau) \right)\ dx + \int_{\overline{\Omega}}1\ d\mathcal{E}(\tau) + \int_0^\tau \int_\Omega F\left(\mathbb{D} u\right) + F^* \left(\mathbb{S}\right)\ dxdt + \int_0^\tau \int_{\Gamma} g \overline{\left|u\right|}\ d\Gamma dt \nonumber \\
\leq& \int_\Omega \frac{1}{2} \frac{\left|\left(\rho u \right)_{0}\right|^2}{\rho_{0}} + P \left(\rho_0 \right)\ dx + \int_0^\tau \int_\Omega \rho f \cdot u\ dxdt \label{energyineq}
\end{align}

is satisfied for almost all $\tau \in [0,T]$.
\end{definition}

We briefly explain the connection between this definition and the classical formulation \eqref{classiccont}--\eqref{coulombbc} of the problem. In the concept of dissipative solutions, the Cauchy stress $\mathbb{T}$ from the classical formulation is replaced by the effective Cauchy stress $\mathbb{T}_{\text{eff}}:= \mathbb{S}_{\text{eff}} - p \text{Id}$, where $\mathbb{S}_{\text{eff}}:= \mathbb{S} - \mathcal{R}$ is called the effective viscous stress with the turbulent component $\mathcal{R}$, which is called the Reynolds stress. Since $\mathcal{R} = 0$ in the classical formulation, we can think of $\mathcal{R}$ as an error term. After this modification, the dissipative momentum inequality \eqref{momineq} can be deduced from the classical momentum equation \eqref{classicmom} as follows: The latter identity is multiplied by an arbitrary test function $\phi \in \mathcal{D}([0,T)\times \overline{\Omega})$ satisfying $\phi \cdot \text{n}|_\Gamma = 0$ and integrated by parts to obtain the relation
\begin{align}
-\int_0^T \int_\Omega \rho u \cdot \partial_t \phi \ dxdt - \int_\Omega \left(\rho u\right)_0\phi(0)\ dx = &\int_0^T \int_\Omega \left(\rho u \otimes u \right): \nabla \phi - \mathbb{T}_{\text{eff}}: \nabla \phi + \rho f \cdot \phi \ dxdt \nonumber \\
&+ \int_0^T \int_\Gamma \left( \mathbb{T}_{\text{eff}} \text{n} \right)_\tau \cdot \phi \ d\Gamma dt. \nonumber
\end{align}

To this identity we then apply the Coulomb friction law boundary condition \eqref{boundcond1}, \eqref{coulombbc} - again with $\mathbb{T}$ replaced by $\mathbb{T}_{\text{eff}}$ - and infer that
\begin{align}
-\int_0^T \int_\Omega \rho u \cdot \partial_t \phi \ dxdt - \int_\Omega \left(\rho u\right)_0\phi(0)\ dx \geq &\int_0^T \int_\Omega \left(\rho u \otimes u \right): \nabla \phi - \mathbb{S}: \nabla \phi + p(\rho)\nabla \cdot \phi + \rho f \cdot \phi \ dxdt \nonumber \\
&+ \int_0^T \int_{\Omega} \nabla \phi :\ d\mathcal{R}(t)dt + \int_0^T \int_\Gamma g|u| - g \left| \phi + u \right| \ d\Gamma dt, \nonumber
\end{align}

where the Reynolds stress $\mathcal{R}$ is interpreted as a tensor-valued measure. Finally, we replace the boundary term $|\phi + u|$ on the right-hand side of this inequality by a function $\overline{|\phi + u|}$ and arrive at the desired dissipative momentum inequality \eqref{momineq}. The reason for this last step is of purely mathematical nature: In the Galerkin limit in the approximation method in the proof of our main result we are not able to show strong convergence of the velocity field $u_n$, see Section \ref{gallim} below. In particular, we are not able to relate the weak limit $\overline{|\phi + u|}$ of the function $|\phi + u_n|$ to the limit velocity $u$. Instead, we identify $\overline{|\phi + u|}$ in the sense of Young measures, which is expressed through the relation \eqref{boundaryid} in our definition of dissipative solutions. The question if the momentum inequality \eqref{momineq} also holds true for $|\phi + u|$ instead of $\overline{|\phi + u|}$ remains an open problem.

In the energy inequality \eqref{energyineq} in our dissipative formulation, the quantity $\mathcal{E}$, which is referred to as the energy dissipation defect, can be understood as another error term, connected to the error term $\mathcal{R}$ via the compatibility condition \eqref{compcon}. The inequality between the left-hand side and the right-hand side of \eqref{energyineq} constitutes the actual energy inequality of our problem. It is deduced from the classical formulation \eqref{classiccont}--\eqref{coulombbc} - with $\mathbb{T}$ replaced by $\mathbb{T}_{\text{eff}}$ - by testing the momentum equation \eqref{classicmom} by $u$, subtracting the continuity equation \eqref{classiccont} tested by $\frac{1}{2}|u|^2$ and exploiting the Coulomb friction law boundary condition \eqref{boundcond1}, \eqref{coulombbc}. The second inequality in \eqref{energyineq} instead represents a slightly stronger condition, incorporating the (unidentified) boundary function $\overline{|u|} = \overline{|0+u|}$. The latter inequality is important for our considerations, because it guarantees that for any dissipative solution which is in addition smooth it holds that $\mathcal{E}=0$, $\mathcal{R}=0$ and $\mathbb{S}\in \partial F(\mathbb{D}u)$.

Indeed, suppose that $(\rho, u)$ is a smooth dissipative solution to the problem \eqref{classiccont}--\eqref{coulombbc} in the sense of Definition \ref{disssol}. Then, under exploitation of the continuous dependence (with respect to the $L^q((0,T)\times \Gamma)$-norm) of $\overline{|\phi + u|}$ on $\phi$, it is easy to see that $u$ may be used a test function in the dissipative momentum inequality \eqref{momineq}. In combination with the continuity equation \eqref{conteq}, this leads - similarly as in the derivation of the energy inequality \eqref{energyineq} - to an estimate of the form
\begin{align}
& - \int_\Omega \frac{1}{2} \rho (\tau)|u(\tau)|^2 + P \left(\rho (\tau) \right)\ dx + \int_0^\tau \int_\Omega - \mathbb{S}:\nabla u + \rho f\cdot u\ dxdt + \int_0^\tau \int_\Gamma g |u| - g \overline{|u+u|} \ d\Gamma dt \nonumber \\
\leq& -\int_\Omega \frac{1}{2} \rho_0|u_0|^2 + P \left(\rho_0 \right)\ dx - \int_0^\tau \int_\Omega \rho f \cdot u\ dxdt - \int_0^\tau \int_{\Omega} \nabla u : \ d\mathcal{R}(t)dt \label{estimatebound}
\end{align}

for almost all $\tau \in [0,T]$. Here, we exploit the identification \eqref{boundaryid} of the boundary term $\overline{|u+u|}$ to see that
\begin{align}
\overline{|u+u|}(t,y) = \int_{\mathbb{R}^3} \left| u(t,y) + z \right|\ d\nu_{(t,y)}(z) \leq \int_{\mathbb{R}^3} \left| u(t,y) \right| + \left| z \right|\ d\nu_{(t,y)}(z) = \left| u(t,y) \right| + \overline{|u|}(t,y) \nonumber
\end{align}

for almost all $(t,y) \in (0,T)\times \Gamma$. This allows us to rewrite the estimate \eqref{estimatebound} as
\begin{align}
& - \int_\Omega \frac{1}{2} \rho (\tau)|u(\tau)|^2 + P \left(\rho (\tau) \right)\ dx + \int_0^\tau \int_\Omega - \mathbb{S}:\nabla u + \rho f\cdot u\ dxdt + \int_0^\tau \int_\Gamma - g \overline{|u|}\ d\Gamma dt \nonumber \\
\leq& -\int_\Omega \frac{1}{2} \frac{\left|\left(\rho u \right)_{0}\right|^2}{\rho_{0}} + P \left(\rho_0 \right)\ dx - \int_0^\tau \int_\Omega \rho f \cdot u\ dxdt - \int_0^\tau \int_{\Omega} \nabla u : \ d\mathcal{R}(t)dt. \nonumber
\end{align}

We proceed by adding this estimate to the second inequality in the energy estimate \eqref{energyineq} and infer that
\begin{align}
\int_{\overline{\Omega}} 1\ d\mathcal{E}(\tau) + \int_0^\tau \int_\Omega F\left(\mathbb{D}u \right) + F^*\left( \mathbb{S} \right) - \mathbb{S}: \nabla u\ dxdt \leq - \int_0^\tau \int_{\Omega} \nabla u : \ d\mathcal{R}(t)dt. \nonumber
\end{align}

Now, thanks to the compatibility condition \eqref{compcon} between $\mathcal{E}$ and $\mathcal{R}$, it suffices, as explained in \cite[Section 4]{afn}, to apply the Gronwall inequality to first infer that $\mathcal{E}=0$ and $\mathcal{R}=0$ and subsequently
\begin{align}
F\left(\mathbb{D}u \right) + F^*\left( \mathbb{S} \right) - \mathbb{S}: \nabla u = 0, \nonumber
\end{align}

i.e.\ $\mathbb{S} \in \partial F(\mathbb{D}u)$ according to the characterization \eqref{fenchineq} of subgradients. We are now in the position to present our main result, which proves the existence of dissipative solutions to the problem \eqref{classiccont}--\eqref{coulombbc} in the sense of Definition \ref{disssol}.
\begin{theorem} \label{mainresult}
Let $T > 0$ and let $\Omega \subset \mathbb{R}^3$ be a bounded domain with boundary $\Gamma = \partial \Omega$ of class $C^{2,\eta}\bigcup C^{0,1}$ for some $\eta>0$. Let the data
\begin{align}
f \in L^\infty((0,T)\times \Omega;\mathbb{R}^3),\quad \quad 0 \leq g \in L^\infty((0,T)\times \Gamma;\mathbb{R}),\quad \quad F:\mathbb{R}^{d \times d}_{\operatorname{sym}} \to [0, \infty),\quad \quad p:[0,\infty) \rightarrow \mathbb{R} \nonumber
\end{align}

satisfy the conditions \eqref{fcond1}, \eqref{fcond2} and \eqref{pcond1}--\eqref{pcond3}. Moreover, consider initial data
\begin{align}
\rho_0 \in L^\gamma(\Omega;\mathbb{R}),\quad \rho_0 \geqslant 0,\quad (\rho u)_0 \in L^{\frac{2 \gamma}{\gamma + 1}}\left(\Omega; \mathbb{R}^3\right),\quad \int_{\Omega}^{}\, \left( \frac{1}{2}\, \frac{\vert (\rho u)_0 \vert^2}{\rho_0} + P(\rho_0) \right)\, dx < \infty, \nonumber
\end{align}

where $P$ denotes the pressure potential defined in \eqref{pcond2} and $\gamma>1$ is chosen according to the estimate \eqref{coercP}. Then the problem \eqref{classiccont}--\eqref{coulombbc} admits a dissipative solution $(\rho, u)$ in the sense of Definition \ref{disssol}.
\end{theorem}

We point out that the reason why we assume the domain boundary $\Gamma$ to be of class $C^{2,\eta}$ in this theorem lies in our application of the classical parabolic regularization procedure of the continuity equation (see e.g.\ \cite[Section 7.7]{ns}) for the construction of the density in Section \ref{approxexist} below. The $C^{0,1}$-regularity of $\Gamma$ instead is required for the application of Korn's inequality (see Lemma \ref{korn} in the Appendix) as well as for the theory of Young measures (see Lemma \ref{youngmeasures} in the Appendix). The remainder of the article is dedicated to the proof of Theorem \ref{mainresult}.

\section{Approximate system}\label{sub1}

The proof of Theorem \ref{mainresult} is based on a three level approximation scheme: We first consider an approximate problem consisting of three different levels, each of which is associated to a parameter $\delta, \epsilon >0$, $n \in \mathbb{N}$. This approximate problem is chosen in such a way that it can be solved via classical methods. After solving the approximate problem, we pass to the limit in each approximation level separately and finally recover a solution to the original problem. In the following, we present the full approximate problem and subsequently give a short explanation of the purpose of each individual approximation level.

We choose a sequence $(V_n)_{n \in \mathbb{N}}$ of $n$-dimensional vector spaces 
\begin{align}
V_n \subset C^{2,\eta}(\overline{\Omega};\mathbb{R}^3) \subset L^2(\Omega;\mathbb{R}^3) \label{galerkinspace}
\end{align}

of the dimension $n$, such that
\begin{align}
\bigcup_{n=1}^\infty V_n \quad \text{is dense in} \quad W_{\text{n}}^{2,p}(\Omega) \quad \text{for some} \quad 3 < p < \infty. \label{galerkindensity}
\end{align}

Then our full approximate problem, containing all approximation levels, consists of finding a solution
\begin{align}
\left(\rho_\delta, u_\delta \right) \in \left\lbrace \psi \in C\left([0,T]; C^{2,\eta} \left( \overline{\Omega} \right) \right) \bigcap C^1\left([0,T]; C^{0,\eta} \left( \overline{\Omega} \right) \right):\ \left. \nabla \psi \cdot \text{n} \right|_{\Gamma} = 0 \right\rbrace \times C\left([0,T];V_n \right) \label{approxsol}
\end{align}

to the regularized continuity equation
\begin{align}\label{c22}
\partial_t \rho_\delta + \nabla \cdot (\rho_\delta u_\delta) = \epsilon \Delta \rho_\delta \quad \text{in } [0,T] \times \Omega
\end{align}
 
and the regularized momentum equation
\begin{align}
\int_\Omega \partial_t \left(\rho_\delta u_\delta \right) \cdot \phi \ dx =& \int _\Omega \left( \rho_\delta u_\delta \otimes u_\delta \right): \nabla \phi- \partial F_\delta \left( \mathbb{D} u_\delta \right): \nabla \phi + p \left(\rho_\delta \right) \nabla \cdot \phi \nonumber \\
& + \rho_\delta f\cdot \phi - \epsilon \left( \nabla u_\delta \nabla \rho_\delta \right) \cdot \phi\ dx - \int_{\Gamma} g \partial j_\delta \left(u_\delta \right) \cdot \phi\ d\Gamma \quad \quad \text{in } [0,T] \label{m22}
\end{align}

for any $\phi \in C([0,T]; V_n)$ which in addition satisfies the initial conditions
\begin{align}
\rho_\delta (0,\cdot) = \rho_0 \quad \text{in } \Omega,\quad \quad \left( \rho_\delta u_\delta \right) (0,\cdot) = \left(\rho u\right)_0 \quad \text{in } \Omega. \label{initialcondreg}
\end{align}

In the momentum equation \eqref{m22} of this approximate problem, $F_\delta$ constitutes a regularization of the potential $F$, defined by
\begin{align}
F_\delta (D) = \int_{\mathbb{R}^{3\times 3}_{\operatorname{sym}}}^{}\, \xi_\delta (D - \mathbb{Z} ) F(\mathbb{Z}) d \mathbb{Z} - \int_{\mathbb{R}^{3\times 3}_{\operatorname{sym}}}^{}\ \xi_\delta( \mathbb{Z} ) F(\mathbb{Z}) d \mathbb{Z} \quad \quad \forall D \in \mathbb{R}^{3\times 3}_{\operatorname{sym}} \label{regularpotential}
\end{align}

for a radially symmetric and non-increasing mollifier $\xi_\delta$ . We notice that, due to the smoothness of $F_\delta$, the quantity $\partial F_\delta$ constitutes the classical gradient of $F_\delta$. We further point out that, as the original potential $F$, the regularized potential $F_\delta$ is convex, non-negative and satisfies the relations
\begin{align}\label{Fdelta}
F_\delta (0) = 0,\quad \quad F_\delta(D) \geqslant \mu \Big \vert D - \frac{1}{3} \operatorname{tr} \left[ D\right] \operatorname{Id} \Big \vert^q \quad \forall D \in \mathbb{R}^{3\times 3}_{\operatorname{sym}} \quad \text{with} \quad \vert D \vert >1
\end{align}
for a constant $\mu > 0$ and the same exponent $q>1$ as in \eqref{fprop}. The quantity $j_\delta$ in the momentum equation \eqref{m22} denotes a convex regularization of the absolute value function defined by
\begin{align}
j_\delta (v) := \left\{
\begin{matrix}
|v| & \text{if } |v| > \delta, \\
\frac{|v|^2}{2\delta} + \frac{\delta}{2} & \text{if } |v| \leq \delta.
\end{matrix}\right. \label{J}
\end{align}

This function, which has the regularity
\begin{align}
j_\delta \in C^1(\mathbb{R}^3)\bigcap C^{1,1}_{\operatorname{loc}}(\mathbb{R}^3), \label{jlip}
\end{align}

satisfies the conditions
\begin{align}\nonumber
j_\delta(0) =& 0,\\
\nonumber
\nabla j_\delta(v) \cdot v \geq& 0,\\
\label{jj2} \vert \partial j_\delta(v) \vert \leq& 1,\\
\label{jjj} \vert j_\delta (v) - \vert v \vert \vert \leq& \delta
\end{align}

for all $v \in \mathbb{R}^3$; the quantity $\partial j_\delta$ denotes the gradient of $j_\delta$. Finally, for the initial conditions \eqref{initialcondreg}, the initial data $\rho_0, (\rho u)_0$ is assumed to satisfy the additional conditions
\begin{align}
\rho _0 \in C^{2,\eta}\left(\overline{\Omega}\right),\quad \quad \frac{1}{n} \leq \rho_0 \leq n,\quad \quad \left. \nabla \rho_0 \cdot \operatorname{n} \right|_{\Gamma} = 0, \quad \quad (\rho u)_0 \in C^2 \left( \overline{\Omega} \right). \label{-342}
\end{align}

After having introduced the full approximate problem \eqref{approxsol}--\eqref{initialcondreg}, we are now in the position to briefly explain the ideas behind the approximation scheme. The general structure of our proof closely mimics the one of the proof of the existence of dissipative solutions to general models of compressible viscous flow with a non-linear viscosity tensor given in \cite{afn}. It is that article, from which we adopt the idea for the regularization \eqref{regularpotential} of the potential $F$ on the $\delta$-level. This technique allows us to replace the viscous stress $\mathbb{S}$, which in the original problem constitutes a subgradient of $F$ at $\mathbb{D}u$, by the classical gradient of $F_\delta $ at $\mathbb{D}u$. In this way, we are able to establish the correct relation between $\mathbb{S}$ and $F$ and, in particular, to deduce a meaningful energy inequality. However, the $\delta$-level in our approximation scheme also incorporates an additional approximation technique, namely the replacement of the boundary integrals in the original momentum inequality \eqref{momineq} by the boundary integral
\begin{align}
\int_{\Gamma} g \partial j_\delta \left(u_\delta \right) \cdot \phi\ d\Gamma. \nonumber
\end{align}

After solving the approximate momentum equation, we may then exploit the convexity of $j_\delta$ to turn it into an inequality in which the Coulomb friction law boundary condition is included in the same way as in our dissipative formulation. This procedure is known, for example, from the proof of the existence of weak solutions to the incompressible Navier-Stokes equations with the Coulomb friction law boundary condition, cf.\ \cite{bmt}. %As a minor addition, we further include the quantity $\delta \nabla u_\delta : \nabla \phi$ in the momentum equation \eqref{m22} on the $\delta$-level: Due to its non-linearity, the quantity $\partial F_\delta$ alone does not suffice to deduce the uniform $L^2(0,T;V_n)$-bounds required for the fixed point argument in the construction of the solution to the approximate momentum equation. The quantity $\delta \nabla u_\delta : \nabla \phi$ fixes this problem, cf.\ \eqref{25} below.

On the $\epsilon$-level of our approximation scheme we employ the classical parabolic regularization of the continuity equation (cf.\ \cite[Section 7.6]{ns}), which shows itself in form of the additional Laplacian on the right-hand side of the approximate continuity equation \eqref{c22}. This procedure guarantees non-negativity of the density in the approximate and hence also in the final system. The quantity $\epsilon \nabla u_\delta \nabla \rho_\delta$ in the approximate momentum equation \eqref{m22} serves as a compensation for this technique, making sure that the energy inequality is preserved under this modification. Finally, the classical Galerkin method carried out on the $n$-level has the purpose of allowing us to solve the momentum equation in the approximate problem. We point out that the density of our Galerkin spaces in $W^{2,p}_{\text{n}}(\Omega)$ in \eqref{galerkindensity} is required to make sure that general test functions $\phi \in \mathcal{D}([0,T)\times \overline{\Omega})$ with $\phi \cdot \text{n}|_{\Gamma} = 0$ can be approximated by Galerkin functions in the norm of $W^{1,\infty}((0,T)\times \Omega)$, cf.\ the convergence \eqref{testconv} below.
%\textcolor{blue}{It seems to me that here we have to add the additional term $\delta \nabla_x u: \nabla_x \phi$. The reason for this is the following: When we linearize the momentum equation below (see \eqref{linmomeq}) we have to replace $\partial F_\delta (\mathbb{D}_x u)$ by $\partial F_\delta (\mathbb{D}_x w)$, since $F_\delta$ is nonlinear. Then, when we later have to show bounds of fixed points of $s\mathbb{T}$ uniformly in $s$ (c.f.\ \eqref{25}), the $\partial F_\delta$-term has the additional factor $s$ and so we don't get a bound uniform with respect to $s$. However, we can still get the bound in the usual way if we add the linear term $\delta \nabla_x u: \nabla_x \phi$}

\section{Existence of the approximate solutions} \label{approxexist}

For the proof of the existence of a solution to our approximate system we can follow almost exactly the corresponding proof in the setting of inflow-outflow boundary conditions and a linear pressure in \cite[Section 4.2]{b}. The only differences in our case lie in the additional boundary integral in the momentum inequality and the non-linearity of the pressure. These, however, do not cause any additional difficulties thanks to the (local) Lipschitz-continuity of the pressure (see Remark \ref{pprops}) and the quantity $\partial j_\delta$ (see the inclusion \eqref{jlip}). Nonetheless, for the sake of self-containedness, we sketch the proof in our setting in the present section. Since the proof in \cite[Section 4.2]{b}, in turn, closely resembles the well-known proof of the existence of approximate weak solutions to the compressible Navier-Stokes equations we further refer, for some additional details, to \cite[Section 7.7]{ns}.

As a first step, we take a look at the continuity equation on its own: For an arbitrary but fixed function $v \in C([0,T];V_n)$ and some initial function $\rho_0 \in C^{2,\eta}(\Omega)$, satisfying
\begin{align}
\frac{1}{n} \leq \rho_0 \leq n \quad \text{in } \Omega, \quad \quad \nabla \rho_0 \cdot \text{n} = 0 \quad \text{on } \Gamma, \nonumber
\end{align}

we consider the problem
\begin{align}
\partial_t \rho + \nabla \cdot \left( \rho v \right) = \epsilon \Delta \rho \quad \text{in } [0,T] \times \Omega, \quad \quad \rho(0, \cdot) = \rho_0(\cdot)\quad \text{in } \Omega. \label{-228}
\end{align}

The classical theory for the parabolic regularization of the continuity equation (see \cite[Lemma 3.1, Theorem 10.22, Theorem 10.23]{fn}, \cite[Proposition 7.39]{ns}) guarantees the existence of a unique solution
\begin{align}
\rho = \rho(v) \in \left\lbrace \psi \in C\left([0,T]; C^{2,\eta} \left( \overline{\Omega} \right) \right) \bigcap C^1\left([0,T]; C^{0,\eta} \left( \overline{\Omega} \right) \right):\ \left. \nabla \psi \cdot \text{n} \right|_{\Gamma} = 0 \right\rbrace \nonumber
\end{align}

to this problem. Additionally, for any bounded set $B \subset C([0,T];V_n)$, this solution satisfies the estimates
\begin{align}
\frac{1}{n} \exp \left(-\left\| v \right\|_{L^1(0,T;W^{1,\infty}(\Omega))} \right) \leq \rho(t,\cdot) \leq& n \exp \left(\left\| v \right\|_{L^1(0,T;W^{1,\infty}(\Omega))} \right) \quad \text{in } \overline{\Omega} \label{-233} \\
\left\| \rho(v) \right\|_{C \left([0,T];C^{2,\eta}\left(\overline{\Omega}\right)\right)} + \left\| \rho(v) \right\|_{C^1 \left([0,T];C^{0,\eta}\left(\overline{\Omega}\right)\right)} \leq& c(B), \label{-330} \\
\left\| \rho \left(v^1\right) - \rho \left(v^2\right) \right\|_{C([0,T];L^2(\Omega))} \leq& c\left(B \right) \left\| v^1 - v^2 \right\|_{C([0,T];W^{1,\infty}(\Omega))} \label{-331}
\end{align}

for all $t \in [0,T]$, all $v,v^1,v^2 \in B$ and a constant $c(B)>0$ depending on $B$ but not on the specific choices of $v,v^1,v^2$. Next, we turn to the momentum equation: We search for $u \in C([0,T];V_n)$ satisfying
\begin{align}
\int_\Omega \left( \rho u \right)(t, \cdot) \cdot \phi \ dx - \int_\Omega \left( \rho u \right)_0 \cdot \phi \ dx =& \int_0^t\int _\Omega \left( \rho u \otimes u \right): \nabla \phi- \partial F_\delta \left( \mathbb{D} u \right): \nabla \phi + p \left(\rho \right) \nabla \cdot \phi \nonumber \\
& + \rho f\cdot \phi - \epsilon \left( \nabla u \nabla \rho \right) \cdot \phi\ dx - \int_{\Gamma} g \partial j_\delta(u) \cdot \phi\ d\Gamma dt \label{linmomeq}
\end{align}

for any $t \in [0,T]$ and any $\phi \in V_n$, where $\rho = \rho(u)$ denotes the solution to the associated continuity equation \eqref{-228} with $v=u$. In order to express this equation as a fixed point problem, we introduce the operators
\begin{align}
\mathcal{M}_{\rho (v) (t)}: V_n \rightarrow V_n^*,\quad &\left\langle \mathcal{M}_{\rho (v)(t)}w,\phi \right\rangle_{V_n^* \times V_n} := \int_\Omega \rho (v)(t) w \cdot \phi\ dx, \nonumber \\
\left\langle \mathcal{N}\left(\rho, u \right), \phi \right\rangle_{V_n^*\times V_n} :=& \int _\Omega \left( \rho u \otimes u \right): \nabla \phi - \partial F_\delta \left( \mathbb{D}u \right):\nabla \phi + p \left(\rho \right) \nabla \cdot \phi \nonumber \\
&+ \rho f\cdot \phi - \epsilon \left( \nabla u \nabla \rho \right) \cdot \phi\ dx - \int_{\Gamma} g \partial j_\delta (u) \cdot \phi\ d\Gamma, \nonumber
\end{align}

as well as the notation
\begin{align}
\left\langle \left( \rho u \right)_0^*, \phi \right\rangle_{V_n^*\times V_n} :=& \int _\Omega \left(\rho u\right)_0 \cdot \phi \ dx. \nonumber
\end{align}

Due to the bound \eqref{-233} of $\rho(v)$ away from zero, the operator $\mathcal{M}_{\rho (v)(t)}$ is invertible with an inverse satisfying the estimates
\begin{align}
\left\| \mathcal{M}^{-1}_{\rho (v)(t)} \right\|_{\mathcal{L}(V_n^*,V_n)} &\leq c\left(B\right), \label{18} \\
\left\| \mathcal{M}^{-1}_{\rho (v^1)(t)} - \mathcal{M}^{-1}_{\rho (v^2)(t)} \right\|_{\mathcal{L}(V_n^*,V_n)} &\leq c\left(B \right) \left\| \rho\left(v^1\right)(t) - \rho \left(v^2\right)(t) \right\|_{L^1(\Omega)} \label{17} \\
\left\| \mathcal{M}_{\rho (v)(t)}^{-1} \mathcal{M}_{\partial_ t\rho (w)(t)}\mathcal{M}_{\rho (v)(t)}^{-1} \right\|_{\mathcal{L}(V_n^*,V_n)} &\leq c\left(B\right) \left\| \partial_t \rho (v)(t) \right\|_{L^1(\Omega)} \label{22}
\end{align}

for any arbitrary but fixed bounded subset $B$ of $C([0,T];V_n)$, any $t \in [0,T]$, any $v,v^1,v^2 \in B$ and a constant $c(B) >0$, cf.\ \cite[Section 7.7.1]{ns}. Moreover, the relation
\begin{align}
&\partial_t \left\langle \mathcal{M}^{-1}_{\rho (v)(t)}w(t),\phi \right\rangle_{V_n^* \times V_n} \nonumber \\
=& \left\langle \mathcal{M}_{\rho (v)(t)}^{-1}\mathcal{M}_{\partial_t \rho (v)(t)}\mathcal{M}_{\rho (v)(t)}^{-1}v(t) + \mathcal{M}_{\rho (v)(t)}^{-1} \partial_t w(t),\phi \right\rangle_{V_n^* \times V_n} \quad \text{in } \mathcal{D}'(0,T) \label{timederid}
\end{align}

holds true for any $v,w \in C([0,T];V_n)$ and any $\phi \in V_n$, see again \cite[Section 7.7.1]{ns}. The above operators allow us to express the desired equation \eqref{linmomeq} equivalently in the form
\begin{align}
u(t) = \mathcal{T}[u](t) \label{15}
\end{align}

where the mapping $\mathcal{T}$ is defined as
\begin{align}
\mathcal{T}: B\left(0,K \right) \rightarrow C\left( \left[0,\tilde{T}\right];V_n \right),\quad \quad \mathbb{T}[u](t) := \mathcal{M}_{\rho(u)(t)}^{-1} \left[ \left(\rho u\right)_0^* + \int_{0}^t \mathcal{N}\left(\rho (\tau) , u(\tau) \right)\ d\tau \right], \nonumber
\end{align}

the ball
\begin{align}
B\left(0,K \right) := \left\lbrace v \in C\left( \left[0,\tilde{T}\right];V_n \right): \left\| v \right\|_{L^\infty(0,\tilde{T};V_n)} \leq K \right\rbrace \nonumber
\end{align}

consists of those $v \in C( [0,\tilde{T}];V_n)$ whose norm in $L^\infty(0,\tilde{T};V_n)$ is bounded by $K>0$ and the time $0<\tilde{T}\leq T$ will be fixed later. Hence, in order to solve the problem \eqref{linmomeq} on the interval $[0,\tilde{T}]$, it suffices to prove the existence of a fixed point $u \in B(0,K)$ of $\mathcal{T}$. To this end we show that $\mathcal{T}$ is a contraction from $B(0,K)$ into itself. From the estimates \eqref{-233}--\eqref{-331}, \eqref{18} and \eqref{17} as well as the local Lipschitz-continuity of both $p$ and $\partial j_\delta$ (cf.\ Remark \ref{pprops} and the inclusion \eqref{jlip}) it is easy to see that
\begin{align}
\left\| \mathcal{T}[u](t) \right\|_{V_n} \leq c_1(n)\exp \left(K\tilde{T}\right) \left(\left\| (\rho u)_0^* \right\|_{V_n^*} +c_2(n,K,T)\tilde{T} \right) \label{selfmap}
\end{align}

and
\begin{align}
\left\| \mathcal{T}\left[u_1\right](t) - \mathcal{T}\left[u_2\right](t) \right\|_{V_n} \leq \tilde{T} c_3(n,K,T) \exp \left(2K\tilde{T} \right) \left\| u_1(t) - u_2(t) \right\|_{V_n} \label{contract}
\end{align}

for all $u,u_1,u_2 \in B(0,K)$ and all $t \in [0,\tilde{T}]$, cf.\ \cite[Section 4.2.2]{b}. Choosing $K>0$ sufficiently large such that
\begin{align}
c_1(n)\text{exp}(1) \left\| (\rho u)_0^* \right\|_{V_n^*} \leq \frac{K}{2} \nonumber
\end{align}

and, subsequently, $\tilde{T}>0$ sufficiently small such that
\begin{align}
K\tilde{T} \leq 1,\quad \quad c_1(n)\exp \left(1\right) c_2(n,K,T)\tilde{T} \leq \frac{K}{2},\quad \quad \tilde{T} c_3(n,K,T) \exp \left(2 \right) < 1, \nonumber
\end{align}

we infer from the estimates \eqref{selfmap} and \eqref{contract} that $\mathcal{T}$ is indeed a contraction from $B(0,K)$ into itself. Hence, the Banach fixed point theorem implies the existence of a fixed point $u \in B(0,K)$ of $\mathcal{T}$, which constitutes the desired solution to the equation \eqref{linmomeq} on the interval $[0,\tilde{T}]$. In particular, with the couple $(\rho_\delta, u_\delta):=(\rho(u),u)$ we have found the solution to the approximate continuity equation \eqref{c22} and to the approximate momentum equation \eqref{m22} on $[0,\tilde{T}]$. It remains to check that the interval $[0,\tilde{T}]$ can be extended to the interval $[0,T]$. This is true if we can show that the norm of $u_\delta (t)$ in $V_n$ remains bounded on $[0,\tilde{T}]$ independently of $\tilde{T}$, in which case the above local existence argument may be iterated until we reach the existence time $\tilde{T}=T$. Consequently, we proceed by deriving an energy equation. To this end we subtract the continuity equation \eqref{c22}, tested by $\frac{1}{2}|u_\delta|^2$, from the momentum equation \eqref{m22}, tested by $u_\delta$. The resulting identity is then integrated with respect to time. Under exploitation of the equality
\begin{align}
\int_\Omega p \left( \rho_\delta \right) \nabla \cdot u_\delta\ dx = - \int_\Omega \epsilon P''\left(\rho_\delta \right) \left| \nabla \rho_\delta \right|^2\ dx - \int_\Omega \partial_t P\left(\rho_\delta \right)\ dx, \nonumber
\end{align}

which follows from the relation \eqref{Pid} between the pressure and its potential, this procedure leads to the energy equation
\begin{align}
&\int_\Omega \frac{1}{2} \rho_\delta (\tau)|u_\delta(\tau)|^2 + P \left(\rho_\delta (\tau) \right)\ dx + \int_0^\tau \int_\Omega \partial F_\delta \left(\mathbb{D}u_\delta \right):\mathbb{D} u_\delta  +\delta \left| \nabla u_\delta \right|^2  + \epsilon P'' \left(\rho_\delta \right) \left| \nabla \rho_\delta \right|^2\ dxdt\nonumber \\
&+ \int_0^\tau \int_{\Gamma} g \partial j_\delta \left(u_\delta\right) \cdot u_\delta\ d\Gamma dt =\int_\Omega \frac{1}{2} \frac{\left|\left(\rho u \right)_{0}\right|^2}{\rho_{0}} + P \left(\rho_0 \right)\ dx + \int_0^\tau \int_\Omega \rho_\delta f \cdot u_\delta\ dxdt \label{energy2}
\end{align}

for all $\tau \in [0,\tilde{T}]$. From this identity, the convexity of $F_\delta$ and the lower bounds \eqref{coercP} for the pressure potential and \eqref{Fdelta} for $F_\delta$ we see that
\begin{align}
\left\| \rho_\delta \left| u_\delta \right|^2 \right\|_{L^\infty(0,\tilde{T};L^1(\Omega))} + \left\| \rho_\delta \right\|_{L^\infty(0,\tilde{T};L^\gamma (\Omega))} + \left\| \mathbb{D}u_\delta - \frac{1}{3} \operatorname{tr} \left[\mathbb{D}u_\delta \right] \operatorname{id} \right\|_{L^q((0,\tilde{T})\times \Omega)} \leq c \nonumber
\end{align}

for a constant $c > 0$ independent of $\tilde{T}$ (and $\delta$). Consequently, the Korn inequality \eqref{kornest} in Lemma \ref{korn} in the Appendix and the equivalence of norms on the $n$-dimensional Galerkin space $V_n$ imply that
\begin{align}
\left\|u_\delta \right\|_{L^q(0,\tilde{T};V_n)} \leq c\left\|u_\delta \right\|_{L^q(0,\tilde{T};W^{1,q}(\Omega))} \leq c. \label{ubound}
\end{align}

This, in turn, in combination with the estimates \eqref{-233}, \eqref{-330} for the solution to the parabolic Neumann problem \eqref{-228} with $v=u_\delta$ and the energy inequality \eqref{energy2} shows that
\begin{align}
\left\| u_{\delta} \right\|_{C([0,\tilde{T}];V_n)} + \left\| \rho_\delta \right\|_{C \left([0,\tilde{T}];C^{2,\eta}\left(\overline{\Omega}\right)\right)} + \left\| \rho_\delta \right\|_{C^1 \left([0,\tilde{T}];C^{0,\eta}\left(\overline{\Omega}\right)\right)} + \left\| \frac{1}{\rho_\delta} \right\|_{C([0,\tilde{T}]\times \overline{\Omega})} \leq c. \label{deltabounds}
\end{align}

This proves the desired bound of $u_\delta (t)$ in $V_n$ independent of $\tilde{T}$ and hence, as stated above, shows that we may choose $\tilde{T}=T$. Finally, for the later use in the limit passage with respect to $\delta \rightarrow 0$, we point out that a combination of the identities \eqref{timederid} and \eqref{15} with the estimates \eqref{-233}, \eqref{-330}, \eqref{18}, \eqref{22} and \eqref{deltabounds} results in the bound
\begin{align}
\left\| \partial_t u_\delta \right\|_{L^2(0,\tilde{T};V_n)} \leq c \label{timederest}
\end{align}

for another constant $c>0$ independent of $\delta$. Altogether, we have shown the following result:
\begin{proposition} 
\label{System on delta level}
Let the conditions of Theorem \ref{mainresult} be satisfied, let $\delta, \epsilon >0$ and let $n \in \mathbb{N}$. Assume the initial data $\rho_0$, $(\rho u)_0$ to satisfy the additional regularity conditions \eqref{-342}. Finally, let the regularized potential $F_\delta$ be defined by \eqref{regularpotential} and let the regularization $j_\delta$ of the absolute value function be defined by \eqref{J}. Then there exists a solution
\begin{align}
\left(\rho_\delta, u_\delta \right) \in \left\lbrace \psi \in C\left([0,T]; C^{2,\eta} \left( \overline{\Omega} \right) \right) \bigcap C^1\left([0,T]; C^{0,\eta} \left( \overline{\Omega} \right) \right):\ \left. \nabla \psi \cdot \text{n} \right|_{\Gamma} = 0 \right\rbrace \times C\left([0,T];V_n \right) \nonumber
\end{align}

to the approximate problem \eqref{approxsol}--\eqref{initialcondreg} which furthermore obeys the energy equation \eqref{energy2} and the bounds \eqref{deltabounds}, \eqref{timederest} for $\tilde{T}=T$ and a constant $c>0$ independent of $\delta$.
\end{proposition}

\section{Limit passage with respect to \texorpdfstring{$\delta \rightarrow 0$}{}} \par \label{reglim}

Next, we let $\delta$ tend to zero in order to pass to the limit in the regularizations $F_\delta$ and $j_\delta$ of the potential $F$ and the absolute value function, respectively. Under exploitation of the bounds \eqref{deltabounds}, \eqref{timederest} - which, according to Proposition \ref{System on delta level}, hold true for $\tilde{T}=T$ and independently of $\delta$ - and the Aubin-Lions lemma we find functions
\begin{align}
0 \leq \rho \in& \bigg\lbrace \psi \in C\left([0,T];H^{1}(\Omega)\right) \bigcap C\left([0,T];L^r(\Omega) \right) \bigcap L^2 \left(0,T; H^{2}(\Omega) \right) \bigcap L^\infty\left((0,T)\times \Omega \right): \nonumber \\
&\ \ \partial_t \psi \in L^2 \left((0,T)\times \Omega \right),\ \frac{1}{\psi} \in L^\infty \left((0,T)\times \Omega \right),\ \left. \nabla \psi \cdot \operatorname{n} \right|_{\partial \Omega} = 0 \bigg\rbrace, \nonumber \\
u \in& \left\lbrace \phi \in C\left([0,T];V_n \right):\ \partial_t \phi \in L^2\left(0,T;V_n \right) \right\rbrace \nonumber
\end{align}

for all $1 \leq r < \infty$, such that, possibly after the extraction of subsequences,
\begin{align}
\rho_{\delta} \rightarrow \rho \quad &\text{in } C\left([0,T];H^{1}(\Omega) \right) \ \ \text{and} \ \ C\left([0,T];L^r( \Omega) \right), \ \quad u_{\delta} \rightarrow u \quad &\text{in }& C\left([0,T];V_n \right), \label{deltacon1} \\
\rho_{\delta} \rightharpoonup \rho \quad &\text{in } L^2\left(0,T;H^{2}(\Omega) \right),\quad \quad \quad \quad \quad \quad \quad \quad \quad \quad \quad \quad \ \rho_\delta \buildrel\ast\over\rightharpoonup \rho \quad &\text{in }& L^\infty \left((0,T)\times \Omega \right), \nonumber \\
\partial_t \rho_{\delta} \rightharpoonup \partial_t \rho \quad &\text{in } L^2\left((0,T)\times \Omega \right),\quad \quad \quad \quad \quad \quad \quad \quad \quad \quad \quad \quad \partial_t u_{\delta} \rightharpoonup \partial_t u \quad &\text{in }& L^2\left(0,T;V_n \right) \label{deltacon3}
\end{align}

for all $1 \leq r < \infty$. These convergences allow us to pass to the limit in the continuity equation \eqref{c22} and infer that
\begin{align}\label{c33}
\partial_t \rho + \nabla \cdot (\rho u) = \epsilon \Delta \rho \quad \text{in } [0,T] \times \Omega.
\end{align}

We point out that, as a consequence, $\rho$ and $u$ also satisfy the renormalized continuity equation
\begin{align}
\partial_t \zeta \left(\rho \right) + \nabla \cdot \left( \zeta \left(\rho \right)u \right) + \left[ \zeta' \left(\rho \right)\rho - \zeta \left(\rho \right) \right] \nabla \cdot u - \epsilon \Delta \zeta \left(\rho \right) = - \epsilon \zeta'' \left(\rho \right)\left| \nabla \rho \right|^2 \leq 0 \quad \text{in } [0,T] \times \Omega \label{rc33}
\end{align}

for all convex functions $\zeta \in C^2([0,+\infty))$. In order to pass to the limit in the momentum equation and the energy inequality, we further need to show convergence of the composite functions in the problem. From the continuous dependence of the pressure $p$ on the density, the uniform convergence of $j_\delta$ implied by the estimate \eqref{jjj} and the strong convergences \eqref{deltacon1} of $\rho_\delta$ and $u_\delta$ we see that
\begin{align}
p\left( \rho_\delta \right) \rightarrow p(\rho) \quad &\text{in } C\left([0,T];L^r(\Omega) \right), \quad \quad \quad \quad \quad \quad \ \ P\left( \rho_\delta \right) \rightarrow P(\rho) \quad &\text{in }& C\left([0,T];L^r(\Omega) \right), \label{deltacon4} \\
j_\delta \left(u_\delta \right) \rightarrow j(u) \quad &\text{in } C\left([0,T] \times \Gamma \right), \quad \quad \quad \quad  j_\delta \left(\phi + u_\delta \right) \rightarrow j(\phi + u) \quad &\text{in }& C\left([0,T] \times \Gamma \right) \label{deltacon5}
\end{align}

for all $1 \leq r < \infty$ and all $\phi \in C([0,T];V_n)$. The convergence of $\partial F_\delta (\mathbb{D} u_\delta)$ is slightly more delicate: We first notice that
\begin{align}
F_\delta \left( D \right) \leq G\left( D \right) := \sup_{\delta \in [0,1]} F_\delta (D) \quad \quad \forall D \in \mathbb{R}^{3 \times 3}_{\operatorname{sym}},\ \delta \in [0,1]. \label{supbound}
\end{align}

Since the functions $F_\delta$ are proper convex, the function $G$ is a proper convex function on $\mathbb{R}^{3 \times 3}_{\operatorname{sym}}$ as well. This implies that its conjugate $G^*$ is superlinear. Moreoever, due to the well-known properties of conjugate functions (cf.\ for example \cite[Proposition 2.21 (iii)]{rindler}), the estimate \eqref{supbound} implies that
\begin{align}
G^* \left( D \right) \leq F_\delta^* \left( D \right) \quad \quad \forall D \in \mathbb{R}^{3 \times 3}_{\operatorname{sym}},\ \delta \in [0,1]. \label{conjugatestimate}
\end{align}

Additionally, we may exploit the relation \eqref{fenchineq} between subgradients and conjugate functions for the classical gradient of the smooth function $F_\delta$ to rewrite
\begin{align}
\int_{0}^{T}\,\int_{\Omega}^{}\,\partial F_\delta \left(\mathbb{D} u_\delta \right): \mathbb{D} u_\delta \ dxdt = \int_{0}^{T}\,\int_{\Omega}^{}\, F_\delta \left(\mathbb{D} u_\delta \right) + F_\delta^* \left(\partial F_\delta \left(\mathbb{D} u_\delta \right) \right)\ dxdt \nonumber
\end{align}

in the energy inequality \eqref{energy2}. Consequently, the energy inequality \eqref{energy2} and the estimate \eqref{conjugatestimate} imply that
\begin{align}
\int_{0}^{T}\,\int_{\Omega}^{} G^* \left(\partial F_\delta \left(\mathbb{D} u_\delta \right) \right)\ dxdt \leq \int_{0}^{T}\,\int_{\Omega}^{}\, F_\delta^* \left(\partial F_\delta \left(\mathbb{D} u_\delta \right) \right)\ dxdt \leq c \nonumber
\end{align}

for a constant $c>0$ independent of $\delta$. Due to the superlinearity of $G^*$, this allows us to apply the de la Vallée Poussin criterion for equiintegrability, which implies the existence of a function $\mathbb{S} \in L^1((0,T)\times \Omega; \mathbb{R}^{3\times 3}_{\operatorname{sym}})$ such that
\begin{align}
\partial F_\delta \left(\mathbb{D}u_\delta \right) \rightharpoonup \mathbb{S} \quad \text{in } L^1\left((0,T)\times \Omega \right). \label{convergences}
\end{align}

We want to identify the limit function $\mathbb{S}$ as a subgradient of $F$ at $\mathbb{D}u$. We achieve this goal by following the arguments from \cite[Section 3.1]{woznicki}: On the $\delta$-level, $\partial F_\delta (\mathbb{D}u_\delta)$ constitutes the classical gradient and in particular a subgradient of $F_\delta$ at $\mathbb{D}u_\delta$. Thus, for all non-negative functions $\phi \in \mathcal{D}((0,T)\times \Omega)$ and all matrices $D \in \mathbb{R}^{3 \times 3}_{\operatorname{sym}}$ it holds that
\begin{align}
0 \leq \int_0^T \int_\Omega \left[ F_\delta \left(D \right) - \partial F_\delta (\mathbb{D}u_\delta) : \left(D - \mathbb{D} u_\delta \right) - F_\delta \left(\mathbb{D} u_\delta \right) \right] \phi \ dxdt \label{nonneg}
\end{align}

Further, from the definition \eqref{regularpotential} of the regularized potential $F_\delta$ we know that
\begin{align}
F_\delta (\cdot) \rightarrow F(\cdot) - F(0) = F(\cdot) \quad \text{in } C_{\text{loc}}\left( \mathbb{R}^{3 \times 3}_{\operatorname{sym}} \right). \nonumber
\end{align} %lokal glm. konvergenz da F stetig ist, siehe zb dobrowolski buch

This, in combination with the uniform convergence \eqref{deltacon1} of $\mathbb{D}u_\delta$ and the weak convergence \eqref{convergences}, allows us to pass to the limit on the right-hand side of the inequality \eqref{nonneg}. Due to the arbitrary choice of $0 \leq \phi \in \mathcal{D}((0,T)\times \Omega)$ we then infer that
\begin{align}
0 \leq F(D) - \mathbb{S}: \left( D - \mathbb{D}u \right) - F \left( \mathbb{D}u \right) \quad \quad \forall D \in \mathbb{R}^{3 \times 3}_{\operatorname{sym}}, \nonumber
\end{align} 

i.e.\ the desired inclusion $\mathbb{S} \in \partial F(\mathbb{D}u)$. In particular, by the equivalent characterization \eqref{fenchineq} of subgradients, it holds that
\begin{align}
\mathbb{S} : \mathbb{D} u = F\left(\mathbb{D} u\right) + F^* \left(\mathbb{S}\right). \label{identitys2}
\end{align}

We have now shown all the convergences necessary for the limit passages in both the momentum equation and the energy inequality. In order to carry out these limit passages, we first notice that, thanks to the convexity of $j_\delta$,
\begin{align}
g \partial j_\delta (u_\delta) \cdot \phi = g \partial j_\delta (u_\delta) \cdot \left(\phi + u_\delta - u_\delta \right) \leq gj_\delta \left(\phi + u_\delta \right) - gj_\delta \left( u_\delta \right) \label{convexest}
\end{align}

for all $\phi \in C([0,T];V_n)$. We apply this estimate to the boundary integral in the momentum equation \eqref{m22} and integrate the resulting inequality to see that
\begin{align}
\int_0^T \int_\Omega \partial_t \left(\rho_\delta u_\delta \right) \cdot \phi \ dxdt \geq& \int_0^T \int _\Omega \left( \rho_\delta u_\delta \otimes u_\delta \right): \nabla \phi- \partial F_\delta \left( \mathbb{D} u_\delta \right): \nabla \phi -\delta \nabla u_\delta : \nabla \phi + p \left(\rho_\delta \right) \nabla \cdot \phi \nonumber \\
& + \rho_\delta f\cdot \phi - \epsilon \left( \nabla u_\delta \nabla \rho_\delta \right) \cdot \phi\ dxdt + \int_0^T \int_{\Gamma} g j_\delta \left(u_\delta \right) - gj_\delta \left( \phi + u_\delta \right) \ d\Gamma dt \nonumber
\end{align}

for all $\phi \in C([0,T];V_n)$. The convergences \eqref{deltacon1}, \eqref{deltacon3}, \eqref{deltacon4}, \eqref{deltacon5} and \eqref{convergences} allow us to pass to the limit in this relation and so we infer that
\begin{align}
\int_0^T \int_\Omega \partial_t \left(\rho u \right) \cdot \phi \ dxdt \geq& \int_0^T \int _\Omega \left( \rho u \otimes u \right): \nabla \phi- \mathbb{S}: \nabla \phi + p \left(\rho \right) \nabla \cdot \phi + \rho f\cdot \phi - \epsilon \left( \nabla u \nabla \rho \right) \cdot \phi\ dxdt \nonumber \\
& + \int_0^T \int_{\Gamma} g \left|u \right| - g \left| \phi + u \right| \ d\Gamma dt \label{m33}
\end{align}

for all $\phi \in C([0,T];V_n)$. For the limit passage in the energy inequality we proceed similarly. Applying the estimate \eqref{convexest} for the choice $\phi = -u_\delta$ to the boundary integral in the energy inequality \eqref{energy2}, we obtain the inequality
\begin{align}
&\int_\Omega \frac{1}{2} \rho_\delta (\tau)|u_\delta(\tau)|^2 + P \left(\rho_\delta (\tau) \right)\ dx + \int_0^\tau \int_\Omega \partial F_\delta \left(\mathbb{D}u_\delta \right):\mathbb{D} u_\delta  +\delta \left| \nabla u_\delta \right|^2  + \epsilon P'' \left(\rho_\delta \right) \left| \nabla \rho_\delta \right|^2\ dxdt\nonumber \\
&+ \int_0^\tau \int_{\Gamma} g j_\delta \left(u_\delta\right)\ d\Gamma dt \leq \int_\Omega \frac{1}{2} \frac{\left|\left(\rho u \right)_{0}\right|^2}{\rho_{0}} + P \left(\rho_0 \right)\ dx + \int_0^\tau \int_\Omega \rho_\delta f \cdot u_\delta\ dxdt \nonumber
\end{align}

for all $\tau \in [0,T]$. On the left-hand side of this inequality we drop the non-negative (cf.\ \eqref{convexity}) quantity $\epsilon P''(\rho_\delta)|\nabla \rho_\delta|^2$. In the resulting inequality we pass to the limit under exploitation of the convergences \eqref{deltacon1}, \eqref{deltacon4}, \eqref{deltacon5}, \eqref{convergences} and the identity \eqref{identitys2} and infer that
\begin{align}
&\int_\Omega \frac{1}{2} \rho (\tau)|u(\tau)|^2 + P \left(\rho (\tau) \right)\ dx + \int_0^\tau \int_\Omega F\left(\mathbb{D} u\right) + F^* \left(\mathbb{S}\right)\ dxdt + \int_0^\tau \int_{\Gamma} g  \left|u\right|\ d\Gamma dt \nonumber \\
\leq& \int_\Omega \frac{1}{2} \frac{\left|\left(\rho u \right)_{0}\right|^2}{\rho_{0}} + P \left(\rho_0 \right)\ dx + \int_0^\tau \int_\Omega \rho f \cdot u\ dxdt \label{energyonepslevel}
\end{align}

for all $\tau \in [0,T]$. Finally, the strong convergences \eqref{deltacon1} also allow us to pass to the limit in the lower and upper bounds \eqref{-233} (for $v = u_\delta$) for the density, showing that the limit density satisfies the corresponding bounds
\begin{align}
\frac{1}{n} \exp \left(-\left\| u \right\|_{L^1(0,T;V_n)} \right) \leq \rho(t,\cdot) \leq& n \exp \left(\left\| u \right\|_{L^1(0,T;V_n)} \right) \quad \text{in } \overline{\Omega} \label{lowupboundseps}
\end{align}

for all $t \in [0,T]$. We end this section with the following remark.
\begin{remark}
Since, by the property \eqref{jj2} of $j_\delta$, it holds that $|\partial j_\delta (u_\delta)| \leq 1$, we find an additional quantity $z \in L^\infty \left((0,T)\times \Gamma \right)$ such that
\begin{align}
\partial j_\delta (u_\delta) \buildrel\ast\over\rightharpoonup z \quad \text{in } L^\infty \left((0,T) \times \Gamma \right),\quad \quad \left\| z \right\|_{L^\infty((0,T)\times \Gamma)} \leq 1. \label{zbound}
\end{align}

Consequently, we may also pass to the limit in the momentum equation \eqref{m22} directly. This leads to the identity
\begin{align}
\int_0^T \int_\Omega \partial_t \left(\rho u \right) \cdot \phi \ dxdt =& \int_0^T\int _\Omega \left( \rho u \otimes u \right): \nabla \phi- \mathbb{S}: \nabla \phi + p \left(\rho \right) \nabla \cdot \phi + \rho f\cdot \phi - \epsilon \left( \nabla u \nabla \rho \right) \cdot \phi\ dxdt \nonumber \\
&- \int_0^T\int_{\Gamma} g z \cdot \phi\ d\Gamma dt \label{malt}
\end{align}

for all $\phi \in C([0,T];V_n)$ and - after an application of the fundamental theorem of the calculus of variations - for all $\phi \in L^\infty(0,T;V_n)$. The identity \eqref{malt} is referred to as the alternative momentum equation on the $\epsilon$-level. It will serve as a technical tool in the derivation of the strong convergence \eqref{strongucon} of the velocity field in the limit passage with respect to $\epsilon \rightarrow 0$ below. In particular, since the identity is not part of our definition of dissipative solutions, we do not need to further specify the quantity $z$.
\end{remark}

\section{Limit passage with respect to \texorpdfstring{$\epsilon \rightarrow 0$}{}} \par

In this section, we pass to the limit in the regularization of the continuity equation. As on the $\delta$-level (cf.\ the estimate \eqref{ubound}), we infer from the energy inequality \eqref{energyonepslevel} and the Korn inequality given by Lemma \ref{korn} in the Appendix that
\begin{align}
\left\|u_\epsilon \right\|_{L^q(0,T;V_n)} \leq c\left\|u_\epsilon \right\|_{L^q(0,T;W^{1,q}(\Omega))} \leq c \nonumber
\end{align}

for a constant $c>0$ independent of $\epsilon$. From this bound, the upper and lower bounds \eqref{lowupboundseps} for the density and the energy inequality \eqref{energyonepslevel} we then infer that
\begin{align}
\left\|u_\epsilon \right\|_{C([0,T];V_n)} + \left\| \rho_\epsilon \right\|_{L^\infty ((0,T)\times \Omega)} + \left\| \frac{1}{\rho_\epsilon} \right\|_{L^\infty ((0,T)\times \Omega)} \leq c. \nonumber
\end{align}

Additionally, testing the continuity equation \eqref{c33} by $\rho_\epsilon$, we see that
\begin{align}
\epsilon^\frac{1}{2} \left\| \nabla \rho_\epsilon \right\|_{L^2((0,T)\times \Omega)} \leq c. \nonumber
\end{align}

Under exploitation of these bounds and the continuity equation \eqref{c33} we find functions
\begin{align}
0 \leq \rho \in L^\infty \left((0,T)\times \Omega \right),&\quad \quad u \in L^\infty \left(0,T;V_n \right) \nonumber
\end{align}

such that, possibly after the extraction of subsequences,
\begin{align}
\rho_\epsilon \buildrel\ast\over\rightharpoonup \rho \quad &\text{in } L^\infty \left((0,T)\times \Omega \right), \quad \quad \rho_\epsilon \rightarrow \rho \quad &\text{in }& C_{\text{weak}}\left([0,T];L^2(\Omega) \right) \ \text{and} \ L^2\left(0,T;\left(H^1(\Omega) \right)^* \right), \label{epscon1} \\
u_\epsilon \buildrel\ast\over\rightharpoonup u \quad &\text{in } L^\infty \left(0,T;V_n \right), \quad \quad \epsilon \nabla \rho_\epsilon \rightarrow 0 &\text{in }& L^2 \left((0,T)\times \Omega \right). \label{epscon2} \end{align}

These convergences are sufficient to pass to the limit in the continuity equation \eqref{c33} and infer that the limit functions $\rho$ and $u$ solve the continuity equation
\begin{align}
\partial_t \rho + \nabla \cdot \left(\rho u \right) = 0 \quad \text{in } \mathcal{D}'\left((0,T)\times \Omega \right). \label{c44}
\end{align}

Additionally, the regularization technique by DiPerna and Lions (see for example \cite[Lemma 6.9]{ns}) then shows that $\rho$ and $u$ also satisfy the renormalized continuity equation
\begin{align}
\partial_t \zeta (\rho) + \nabla \cdot \left( \zeta \left(\rho \right)u \right) + \left[ \zeta'\left(\rho \right)\rho - \zeta \left(\rho \right) \right] \nabla \cdot u = 0 \quad \text{in } \mathcal{D}' \left((0,T) \times \mathbb{R}^3 \right) \label{rc44} 
\end{align}

for all
\begin{align}
\zeta \in C^1\left( [0,\infty) \right) \quad \text{such that} \quad \left|\zeta'(r)\right| \leq cr^{\sigma} \quad \forall r \geq 1 \quad \text{for certain } c>0,\ \sigma > -1. \nonumber
\end{align}

For the limit passages in the momentum inequality and the energy inequality, we first notice that the energy inequality \eqref{energyonepslevel}, the superlinearity of $F^*$ and the De La Vallée Poussin criterion for equiintegrability further imply the existence of
\begin{align}
\mathbb{S} \in L^1\left((0,T)\times \Omega; \mathbb{R}^{3 \times 3}_{\text{sym}} \right), \nonumber
\end{align}

such that
\begin{align}
\mathbb{S}_\epsilon \rightharpoonup \mathbb{S} \quad \text{in } L^1 \left((0,T)\times \Omega \right). \label{sconveps}
\end{align}

Moreover, in order to pass to the limit in the pressure and the pressure potential, we show strong $L^p$-convergence of the density. To this end we notice that, due to the bound \eqref{lowupboundseps} of the density away from zero, we may choose $\zeta (\xi) := \xi \ln (\xi)$ in the renormalized continuity equation \eqref{rc33} on the $\epsilon$-level. This choice leads to the inequality
\begin{align}
\int_0^\tau \int_\Omega \rho_\epsilon \operatorname{div} u_\epsilon \ dxdt \leq \int_\Omega \rho_0 \ln \left(\rho_0 \right)\ dx - \int_\Omega \rho_\epsilon (\tau) \ln \left(\rho_\epsilon (\tau) \right)\ dx \label{ineq}
\end{align}

for $\tau \in [0,T]$. After a dominated convergence argument, the same choice is possible in the renormalized continuity equation \eqref{rc44} in the limit, yielding the identity
\begin{align}
\int_0^\tau \int_\Omega \rho \operatorname{div} u \ dxdt = \int_\Omega \rho_0 \ln \left(\rho_0 \right)\ dx - \int_\Omega \rho (\tau) \ln \left(\rho (\tau) \right)\ dx. \label{eq}
\end{align}

Since $\rho_\epsilon$ converges strongly in $L^2(0,T;(H^1(\Omega))^*)$ (cf.\ the convergences \eqref{epscon1}) and $\operatorname{div}u_\epsilon$ converges weakly in $L^2(0,T;H^1(\Omega))$ (cf.\ the convergences \eqref{epscon2}), we further know that
\begin{align}
\rho_\epsilon \operatorname{div} u_\epsilon \buildrel\ast\over\rightharpoonup \rho \operatorname{div} u \quad \text{in } L^\infty \left((0,T) \times \Omega \right). \nonumber
\end{align}

With this knowledge at hand, we may pass to the limit in the inequality \eqref{ineq} and compare the result to the identity \eqref{eq} to see that
\begin{align}
\int_\Omega \rho (\tau) \ln \left(\rho (\tau) \right)\ dx \geq \int_\Omega \overline{\rho \ln (\rho)}(\tau) \ dx \nonumber
\end{align}

for $\tau \in [0,T]$, where $\overline{\rho \ln (\rho)}$ represents a weak $L^1((0,T)\times \Omega)$-limit of $\rho_\epsilon \ln(\rho_\epsilon)$. On the other hand, from the (strict) convexity of the mapping $\xi \mapsto \xi \ln (\xi)$ and the classical relations between convex functions and weak convergence (see \cite[Theorem 10.20]{fn}), we already know that
\begin{align}
\rho \ln (\rho) \leq \overline{\rho \ln (\rho)} \quad \text{in } (0,T)\times \Omega. \nonumber
\end{align}

A combination of the latter two inequalities shows that in fact $\rho \ln (\rho) = \overline{\rho \ln (\rho)}$ and so, again due to the relations between convex functions and weak convergence given by \cite[Theorem 10.20]{fn}, we infer that
\begin{align}
\rho_\epsilon \rightarrow \rho \quad \text{in } L^r \left((0,T)\times \Omega \right) \label{strongdenscon}
\end{align}

for all $1 \leq r < \infty$. In particular, due to the continuous dependence of the pressure on the density $\rho$, we see that
\begin{align}
p\left(\rho_\epsilon \right) \rightarrow p \left(\rho \right) \quad \text{in } L^r\left((0,T)\times \Omega\right),\quad \quad P\left(\rho_\epsilon \right) \rightarrow P \left(\rho \right) \quad \text{in } L^r\left((0,T)\times \Omega\right) \nonumber
\end{align}

for all $1 \leq r < \infty$. Next, for the limit passage in the convective term as well as the boundary terms, we show strong convergence of the velocity field. To this end, we test the alternative momentum equation \eqref{malt} by an arbitrary function $\phi \in L^\infty (0,T;V_n)$. Since $\mathbb{S}_\epsilon$ is bounded in $L^1((0,T)\times \Omega)$ (cf.\ the convergence \eqref{sconveps}) and $|z_\epsilon|\leq 1$ by the bound \eqref{zbound}, this yields the estimate
\begin{align}
&\int_0^T \int_\Omega \partial_t \left( \rho_\epsilon u_\epsilon \right) \cdot \phi \ dxdt \nonumber \\
=& \int_0^T \int _\Omega \left( \rho_\epsilon u_\epsilon \otimes u_\epsilon \right): \nabla \phi- \mathbb{S}_\epsilon :\nabla \phi + p \left(\rho_\epsilon \right) \nabla \cdot \phi + \rho_\epsilon f\cdot \phi - \epsilon \left( \nabla u_\epsilon \nabla \rho_\epsilon \right) \cdot \phi\ dxdt - \int_0^T \int_{\Gamma} g z_\epsilon \cdot \phi\ d\Gamma dt \nonumber \\
\leq& c \left\| \phi \right\|_{L^\infty (0,T;V_n)} \nonumber
\end{align}

for a constant $c>0$ independent of $\epsilon$. Consequently, it holds that
\begin{align}
\left\| \partial_t \left( \rho_\epsilon u_\epsilon \right) \right\|_{(L^\infty (0,T;V_n))^*} \leq c. \label{dualest}
\end{align}

Moreover, by the regularity \eqref{deltacon3} on the $\epsilon$-level, we know that $\partial_t (\rho_\epsilon u_\epsilon) \in L^2((0,T)\times \Omega) \subset L^1(0,T;V_n^*)$. Since the norms of the spaces $L^1(0,T;V_n^*)$ and $(L^\infty (0,T;V_n))^*$ coincide, the bound \eqref{dualest} thus implies that
\begin{align}
\left\| \partial_t \left( \rho_\epsilon u_\epsilon \right) \right\|_{L^1 (0,T;V_n^*)} \leq c. \nonumber
\end{align}

We may therefore apply the Aubin-Lions lemma and infer that
\begin{align}
\rho_\epsilon u_\epsilon \rightarrow \rho u \quad \text{in } L^2 \left(0,T;V_n^* \right). \nonumber
\end{align}

As a consequence we see that
\begin{align}
&\left| \int_0^T \int_\Omega \rho_\epsilon \left| u_\epsilon \right|^2\ dxdt - \int_0^T \int_\Omega \rho \left| u \right|^2\ dxdt \right| \nonumber \\
\leq& \left\| \rho_\epsilon u_\epsilon - \rho u \right\|_{L^2(0,T;V_n^*)} \left\| u_\epsilon \right\|_{L^2(0,T;V_n)} + \left| \int_0^T \int_\Omega \rho u \cdot \left( u_\epsilon - u \right)\ dxdt \right| \rightarrow 0. \label{normconvergence}
\end{align}

In addition, the strong convergence \eqref{strongdenscon} of the density and its bound \eqref{lowupboundseps} away from zero show that $\sqrt{\rho_\epsilon} u_\epsilon$ converges to $\sqrt{\rho}u$ weakly in $L^2((0,T)\times \Omega)$. This, in combination with the convergence \eqref{normconvergence} implies that
\begin{align}
\sqrt{\rho_\epsilon} u_\epsilon \rightarrow \sqrt{\rho} u \quad \text{in } L^2 \left((0,T)\times \Omega \right). \nonumber
\end{align}

Due to the strong convergence \eqref{strongdenscon} of the density, its bound \eqref{lowupboundseps} away from zero and the equivalence of norms on the finite dimensional space $V_n$ this shows the desired strong convergence
\begin{align}
u_\epsilon \rightarrow u \quad \text{in } L^r \left(0,T;V_n \right) \label{strongucon}
\end{align}

for all $1 \leq r < \infty$. In particular, we conclude the convergences
\begin{align}
\left| u_\epsilon \right| \rightarrow \left| u \right| \quad \text{in } L^r \left((0,T)\times \Gamma \right),\quad \quad \left| u_\epsilon + \phi \right| \rightarrow \left| u +\phi \right| \quad \text{in } L^r \left((0,T)\times \Gamma \right) \label{epscon4}
\end{align}

for all $1 \leq r < \infty$ and all $\phi \in \mathcal{D}([0,T);V_n)$. Now, making use of the convergences \eqref{epscon2}, \eqref{sconveps}, \eqref{strongdenscon}, \eqref{strongucon} and \eqref{epscon4}, we may pass to the limit in the momentum inequality \eqref{m33} and obtain the corresponding limit relation
\begin{align}
-\int_0^T \int_\Omega \rho u \cdot \partial_t \phi \ dxdt - \int_\Omega \left( \rho u\right)_0 \phi(0)\ dx \geq& \int_0^T \int _\Omega \left( \rho u \otimes u \right): \nabla \phi- \mathbb{S}: \nabla \phi + p \left(\rho \right) \nabla \cdot \phi + \rho f\cdot \phi \ dxdt \nonumber \\
& + \int_0^T \int_{\Gamma} g \left|u \right| - g \left| \phi + u \right| \ d\Gamma dt \label{m44}
\end{align}

for all $\phi \in \mathcal{D}([0,T);V_n)$. Finally, exploiting the same convergences and in addition the weak lower semicontinuity of both $F$ and $F^*$ (cf.\ Remark \ref{fprops}), we may also pass to the limit in the energy inequality \eqref{energyonepslevel} and infer that
\begin{align}
&\int_\Omega \frac{1}{2} \rho (\tau)|u(\tau)|^2 + P \left(\rho (\tau) \right)\ dx + \int_0^\tau \int_\Omega F\left(\mathbb{D} u\right) + F^* \left(\mathbb{S}\right)\ dxdt + \int_0^\tau \int_{\Gamma} g  \left|u\right|\ d\Gamma dt \nonumber \\
\leq& \int_\Omega \frac{1}{2} \frac{\left|\left(\rho u \right)_{0}\right|^2}{\rho_{0}} + P \left(\rho_0 \right)\ dx + \int_0^\tau \int_\Omega \rho f \cdot u\ dxdt \label{energyonnlevel}
\end{align}

for all $\tau \in [0,T]$.

\section{Limit passage with respect to \texorpdfstring{$n \rightarrow \infty$}{}} \par \label{gallim}

Finally, it remains to pass to the limit in the Galerkin approximation. In order to further return from the regularized initial data from the approximate problem to the initial data from Theorem \ref{mainresult}, we consider a sequence $(\rho_{0,n},(\rho u)_{0,n})$ of initial data satisfying the conditions \eqref{-342} for any fixed $n \in \mathbb{N}$ and such that
\begin{align}
\rho_{0,n} \rightarrow \rho_0 \quad \text{in } L^\gamma (\Omega),\quad \quad \frac{\left|(\rho u)_{0,n} \right|^2}{\rho_{0,n}} \rightarrow \frac{\left|(\rho u)_{0} \right|^2}{\rho_{0}} \quad \text{in } L^1(\Omega) \label{inconv}
\end{align}

for $n \rightarrow \infty$, cf.\ \cite[Section 4]{fnp}. From the energy inequality and the Korn inequality given by Lemma \ref{korn} in the Appendix, we see that
\begin{align}
\left\| \sqrt{\rho_n}u_n \right\|_{L^\infty(0,T;L^2(\Omega))} + \left\| P \left( \rho_n \right) \right\|_{L^\infty (0,T;L^1 (\Omega))} + \left\| u_n \right\|_{L^q (0,T;W^{1,q}(\Omega))} \leq c \label{nbounds}
\end{align}

for a constant $c>0$ independent of $n$. From this, the continuity equation \eqref{c44} on the $n$-level and the variant of the Div-Curl Lemma given by Lemma \ref{divcurl} in the Appendix, we infer the existence of functions
\begin{align}
\rho \in C_{\operatorname{weak}}\left([0,T];L^\gamma (\Omega) \right),\quad \quad u \in L^q \left(0,T;W^{1,q}(\Omega) \right) \label{reglimit}
\end{align}

such that, possibly after the extraction of subsequences,
\begin{align}
\rho_n \rightarrow \rho \quad \text{in } C_{\operatorname{weak}}\left([0,T];L^\gamma (\Omega) \right),\quad \quad u_n \rightharpoonup u \quad \text{in }& L^q \left(0,T;W_{\operatorname{n}}^{1,q}(\Omega) \right), \label{uconvn} \\ \rho_n u_n \buildrel\ast\over\rightharpoonup \rho u \quad \text{in } L^\infty \left(0,T;L^\frac{2\gamma}{\gamma +1}(\Omega) \right).& \label{momconvn1}
\end{align}

These convergences allow us to pass to the limit in the continuity equation \eqref{c44} on the $n$-level and infer that the limit functions $\rho$ and $u$ satisfy the continuity equation
\begin{align}
\partial_t \rho + \nabla \cdot \left(\rho u \right) = 0 \quad \text{in } \mathcal{D}'\left((0,T) \times \Omega \right) \label{c55}
\end{align}

as well. For the limit passage in the momentum inequality and the energy inequality, we infer from the energy inequality \eqref{energyonnlevel} and the bounds \eqref{nbounds} the existence of
\begin{align}
\mathbb{S} \in L^1\left((0,T)\times \Omega; \mathbb{R}^{3 \times 3}_{\text{sym}} \right), \quad \overline{ p\left(\rho \right)},\overline{ P\left(\rho \right)}, \overline{\rho |u|^2} \in L^\infty \left(0,T;\mathcal{M}\left( \overline{\Omega} \right) \right),\quad \overline{\rho u \otimes u} \in L^\infty \left(0,T;\mathcal{M}\left( \overline{\Omega};\mathbb{R}^{3 \times 3}_{\text{sym}} \right) \right) \label{sreg}
\end{align}

such that
\begin{align}
\mathbb{S}_\epsilon \rightharpoonup \mathbb{S} \quad &\text{in } L^1 \left((0,T)\times \Omega \right),\quad \quad \quad \quad \quad p\left(\rho_n \right) \buildrel\ast\over\rightharpoonup \overline{p(\rho)} \quad &\text{in }& L^\infty \left(0,T;\mathcal{M}\left( \overline{\Omega} \right) \right), \label{sconvn} \\
P\left(\rho_n \right) \buildrel\ast\over\rightharpoonup \overline{P(\rho)} \quad &\text{in } L^\infty \left(0,T;\mathcal{M}\left( \overline{\Omega} \right) \right), \quad \quad \rho_n \left|u_n\right|^2 \buildrel\ast\over\rightharpoonup \overline{\rho|u|^2} \quad &\text{in }& L^\infty \left(0,T;\mathcal{M}\left( \overline{\Omega} \right) \right), \label{Pconvn} \\
\rho_n u_n \otimes u_n \buildrel\ast\over\rightharpoonup \overline{\rho u \otimes u} \quad &\text{in } L^\infty \left(0,T;\mathcal{M}\left( \overline{\Omega} \right) \right). \label{momconvn}
\end{align}

Since, at this point, we are not able to identify the limit functions in the convergences \eqref{sconvn}--\eqref{momconvn} anymore, we now introduce the error terms $\mathcal{R}$, $\mathcal{E}$, from our definition of dissipative solutions as
\begin{align}
\mathcal{R} := \left[ \overline{\rho u \otimes u} + \overline{p(\rho)} \operatorname{id} \right] - \left[ \rho u \otimes u - p(\rho) \operatorname{id} \right] &\in L^\infty \left(0,T;\mathcal{M}^+\left( \overline{\Omega};\mathbb{R}^{3 \times 3}_{\text{sym}} \right) \right), \label{errorterms1} \\
\mathcal{E} := \left[\frac{1}{2} \overline{\rho |u|^2} + \overline{P(\rho)} \right] - \left[\frac{1}{2} \rho|u|^2 + P(\rho) \right] \quad &\in L^\infty \left(0,T;\mathcal{M}^+\left( \overline{\Omega} \right) \right). \label{errorterms}
\end{align}

From the convexity of the functions $P - \underline{a}p$ and $\overline{a}p - P$ assumed in \eqref{pcond3} we know that
\begin{align}
\underline{a} \left( \overline{p(\rho)} - p(\rho) \right) \leq \overline{P(\rho)} - P(\rho) \quad \quad \text{and} \quad \quad \overline{P(\rho)} - P(\rho) \leq \overline{a} \left( \overline{p(\rho)} - p(\rho) \right). \nonumber
\end{align}

For the error terms defined in \eqref{errorterms1}, \eqref{errorterms} this implies the existence of constants $0 < \underline{d} \leq \overline{d}< \infty$ such that
\begin{align}
\underline{d}\mathcal{E} \leq \operatorname{tr} \left[\mathcal{R} \right] \leq \overline{d}\mathcal{E}. \label{compconlim}
\end{align}

In order to pass to the limit in the boundary integrals, we notice that, by the weak convergence \eqref{uconvn} of $u_n$, the non-negativity of the function $g \in L^\infty((0,T)\times \Gamma)$ and the weak lower semicontinuity of the $L^1((0,T)\times \Gamma)$-norm,
\begin{align}
\int_0^\tau \int_\Gamma g \left| u \right| d\Gamma dt \leq \liminf_{n \rightarrow \infty} \int_0^\tau \int_\Gamma g \left| u_n \right| d\Gamma dt \label{boundaryest}
\end{align}

for all $\tau \in [0,T]$. Moreover, for any test function $\phi \in \mathcal{D}([0,T]\times \overline{\Omega})$, we find a function $\overline{|\phi + u|} \in L^q((0,T)\times \Gamma)$ such that
\begin{align}
\left| \phi + u_n \right| \rightharpoonup \overline{|\phi + u|} \quad \text{in } L^q \left((0,T) \times \Gamma \right). \label{boundconvn}
\end{align}

Here, while we are not able to identify the limit function $\overline{|\phi + u|}$ as $|\phi + u|$, we may identify it in the sense of Young measures: We apply Lemma \ref{youngmeasures} in the Appendix with the choice
\begin{align}
\Psi \left(t,y;z \right) := \left| \phi \left( t,y\right) + z \right| \quad \quad \forall t \in (0,T),\ y \in \Gamma,\ z \in \mathbb{R}^3 \nonumber
\end{align}

to find a family $\{\nu_{(t,y)} \}_{(t,y) \in (0,T)\times \Gamma}$ of probability measures $\nu_{(t,y)}$ on $\mathbb{R}^3$, independent of $\phi$, such that $\overline{|\phi + u|}$ can be identified in the sense
\begin{align}
\overline{|\phi + u|}(t,y) = \int_{\mathbb{R}^3} |\phi\left(t,y\right) + z|\ d\nu_{(t,y)}(z) \quad \quad \text{for a.a. } (t,y) \in (0,T) \times \Gamma. \label{idbound}
\end{align}

With the convergences \eqref{uconvn}--\eqref{sconvn}, \eqref{momconvn} and \eqref{boundconvn}, the definition \eqref{errorterms1} and the inequality \eqref{boundaryest} at hand, we may now pass to the limit in the momentum inequality \eqref{m44} to infer that
\begin{align}
-\int_0^T \int_\Omega \rho u \cdot \partial_t \phi \ dxdt - \int_\Omega \left(\rho u\right)_0 \phi(0)\ dx \geq& \int_0^T \int _\Omega \left( \rho u \otimes u \right): \nabla \phi- \mathbb{S}: \nabla \phi + p \left(\rho \right) \nabla \cdot \phi + \rho f\cdot \phi \ dxdt \nonumber \\
&+ \int_0^T \int_\Omega \nabla \phi :\ d\mathcal{R}(t)dt + \int_0^T \int_{\Gamma} g \left|u \right| - g \overline{\left| \phi + u \right|} \ d\Gamma dt \label{m55}
\end{align}

for all $\phi \in \mathcal{D}([0,T);V_N)$ with $N \in \mathbb{N}$ fixed. It remains to replace the test functions $\phi \in \mathcal{D}([0,T);V_N)$ in this inequality by the more general test functions $\phi \in \mathcal{D}([0,T)\times \overline{\Omega})$ with $\phi \cdot \text{n}|_{\Gamma} = 0$. To this end we notice that, by our choice \eqref{galerkinspace}, \eqref{galerkindensity} of the Galerkin spaces, for each such $\phi$ we may find a sequence $(\phi_N)_{N \in \mathbb{N}}$ of Galerkin functions $\phi_N \in \mathcal{D}([0,T);V_N)$ such that
\begin{align}
\phi_N \rightarrow \phi \quad \text{in } W^{1,\infty}\left((0,T)\times \Omega \right) \label{testconv}
\end{align}

for $N \rightarrow \infty$. In particular, it holds that
\begin{align}
\left\| \overline{\left|\phi + u\right|} - \overline{\left|\phi_N + u\right|} \right\|_{L^q((0,T)\times \Gamma)} \leq& \liminf_{n \rightarrow \infty} \left\| \left| u_n + \phi \right| - \left| u_n + \phi_N \right| \right\|_{L^q((0,T)\times \Gamma)} \nonumber \\
\leq& \liminf_{n \rightarrow \infty} \left\| \phi - \phi_N \right\|_{L^q((0,T)\times \Gamma)} \rightarrow 0 \label{testconv1}
\end{align}

for $N \rightarrow \infty$. From the convergences \eqref{testconv} and \eqref{testconv1} it follows that the momentum inequality \eqref{m55} indeed holds true for all test functions $\phi \in \mathcal{D}([0,T)\times \overline{\Omega})$ with $\phi \cdot \text{n}|_{\Gamma} = 0$. Finally, under exploitation of the convergences \eqref{inconv}, \eqref{uconvn}, \eqref{sconvn}, \eqref{Pconvn}, the definition \eqref{errorterms}, the inequality \eqref{boundaryest} and the weak lower semicontinuity of both $F$ and $F^*$ (cf.\ Remark \ref{fprops}), we also pass to the limit in the energy inequality \eqref{energyonnlevel} and infer that
\begin{align}
&\int_\Omega \frac{1}{2} \rho (\tau)|u(\tau)|^2 + P \left(\rho (\tau) \right)\ dx + \int_{\overline{\Omega}}1\ d\mathcal{E}(\tau) + \int_0^\tau \int_\Omega F\left(\mathbb{D} u\right) + F^* \left(\mathbb{S}\right)\ dxdt + \int_0^\tau \int_{\Gamma} g \left|u\right|\ d\Gamma dt \nonumber \\
&\int_\Omega \frac{1}{2} \rho (\tau)|u(\tau)|^2 + P \left(\rho (\tau) \right)\ dx + \int_{\overline{\Omega}}1\ d\mathcal{E}(\tau) + \int_0^\tau \int_\Omega F\left(\mathbb{D} u\right) + F^* \left(\mathbb{S}\right)\ dxdt + \int_0^\tau \int_{\Gamma} g \overline{\left|u\right|}\ d\Gamma dt \nonumber \\
\leq& \int_\Omega \frac{1}{2} \frac{\left|\left(\rho u \right)_{0}\right|^2}{\rho_{0}} + P \left(\rho_0 \right)\ dx + \int_0^\tau \int_\Omega \rho f \cdot u\ dxdt \label{energylimit}
\end{align}

for almost all $\tau \in [0,T]$.

\subsection{Proof of the main result}

Summarizing our findings from Section \ref{gallim}, we are now in the position to finish the proof of our main result Theorem \ref{mainresult}. The regularity \eqref{regularity1}--\eqref{additionalterms} of the quantities $\rho, u, \rho u, \mathbb{S}, \mathcal{E}$ and $\mathcal{R}$ follows from the relations \eqref{reglimit}, \eqref{sreg}, \eqref{errorterms1}, \eqref{errorterms}. The compatibility condition \eqref{compcon} between $\mathcal{E}$ and $\mathcal{R}$ was shown in \eqref{compconlim}. The identification \eqref{boundaryid} of the boundary terms $\overline{|\phi + u|}$ in the momentum inequality via the theory of Young measures was achieved in \eqref{idbound}, the continuous dependence of these terms on the test functions $\phi$ follows since
\begin{align}
\left\| \overline{\left|\phi_1 + u\right|} - \overline{\left|\phi_2 + u\right|} \right\|_{L^q((0,T)\times \Gamma)} \leq \liminf_{n \rightarrow \infty} \left\| \left| u_n + \phi_1 \right| - \left| u_n + \phi_2 \right| \right\|_{L^q((0,T)\times \Gamma)} \leq \left\| \phi_1 - \phi_2 \right\|_{L^q((0,T)\times \Gamma)}. \nonumber
\end{align}

The continuity equation \eqref{conteq} was obtained in \eqref{c55}. The momentum inequality \eqref{momineq} was shown in \eqref{m55} first for all finite dimensional test functions $\phi \in \mathcal{D}([0,T);V_N)$, $N \in \mathbb{N}$, and subsequently for the more general test functions $\phi \in \mathcal{D}([0,T)\times \overline{\Omega})$ with $\phi \cdot \text{n}|_{\Gamma} = 0$, cf.\ the convergences \eqref{testconv} and \eqref{testconv1}. Finally, the energy inequality \eqref{energyineq} was shown in \eqref{energylimit}, which concludes the proof of Theorem \ref{mainresult}.

\section{Appendix}

In order to derive uniform $L^q(0,T;W^{1,q}(\Omega))$-estimates for the velocity field from the energy inequality, we make use of the following variant of the Korn inequality \cite[Theorem 10.17]{fn}.

\begin{satz} \label{korn}
Assume $\Omega \subset \mathbb{R}^3$ to be a bounded domain of class $C^{0,1}$ and let $1<q<\infty$, $M,K >0$, $\gamma > 1$. Then there exists a constant $c = c(q,M,K,\gamma)> 0$ with the following property: If $\rho \in L^\gamma (\Omega)$ is a function satisfying
\begin{align}
\rho \geq 0 \quad \text{in } \Omega,\quad \quad M \leq \int_\Omega \rho \ dx,\quad \quad \int_\Omega \rho^\gamma \ dx \leq K. \nonumber
\end{align}

then it holds that
\begin{align}
\left\| u \right\|_{W^{1,q}(\Omega)} \leq c \left( \left\| \mathbb{D}u - \frac{1}{3} \operatorname{tr} \left[\mathbb{D}u \right] \operatorname{id} \right\|_{L^q(\Omega)} + \int_\Omega \rho \ dx + \int_\Omega \rho|u|^2\ dx \right) \label{kornest}
\end{align}

for all $u \in W^{1,q}(\Omega;\mathbb{R}^3)$.
\end{satz}

\textbf{Proof} 

The result is proved in \cite[Theorem 10.17]{fn} for the integrals on the right-hand side of \eqref{kornest} replaced by the integral
\begin{align}
\int_\Omega \rho |u|\ dx. \nonumber
\end{align}

The version \eqref{kornest} of this inequality then follows by the simple estimate
\begin{align}
\int_\Omega \rho |u|\ dx = \int_{\left\lbrace x \in \Omega : |u(x)| \leq 1\right\rbrace} \rho |u|\ dx + \int_{\left\lbrace x \in \Omega : |u(x)| > 1 \right\rbrace } \rho |u|\ dx \leq \int_\Omega \rho \ dx + \int_\Omega \rho |u|^2\ dx. \nonumber
\end{align}

$\hfill \Box$

In order to identify the momentum function in the limit passage in the Galerkin method in Section \ref{gallim}, we exploit the following variant of the Div-Curl lemma, see \cite[Lemma 8.1]{afn}:

\begin{satz}\label{divcurl}
Let $T>0$ and let $\Omega \subset \mathbb{R}^3$ be a bounded domain. Let $\gamma,q,r>1$ and let $(\rho_n)_{n \in \mathbb{N}} \subset L^\gamma((0,T)\times \Omega; \mathbb{R})$ and $(u_n)_{n \in \mathbb{N}} \subset L^q((0,T)\times \Omega; \mathbb{R}^3)$ be two sequences of functions such that
\begin{align}
\rho_n \rightharpoonup \rho \quad \text{in } L^\gamma \left((0,T)\times \Omega \right),\quad \quad u_n \rightharpoonup u \quad \text{in }& L^q\left((0,T)\times \Omega \right), \nonumber \\
\rho_nu_n \rightharpoonup \overline{\rho u} \quad \text{in } L^r\left((0,T) \times \Omega \right) \nonumber
\end{align}

for some functions $\rho \in L^\gamma ((0,T)\times \Omega)$, $u \in L^q((0,T)\times \Omega)$, $\overline{\rho u} \in L^r((0,T)\times \Omega)$. In addition, suppose that
\begin{align}
\partial_t \rho_n + \nabla \cdot \left(\rho_n u_n \right) = 0 \quad \text{in } \mathcal{D}'\left((0,T) \times \Omega \right). \nonumber 
\end{align}

and
\begin{align}
\left\| \nabla u_n \right\|_{L^1((0,T)\times \Omega)} \leq c \nonumber
\end{align}

for a constant $c>0$ independent of $n$. Then it holds that
\begin{align}
\overline{\rho u} = \rho u \quad \text{a.e. in } (0,T) \times \Omega. \nonumber
\end{align}

\end{satz}

\textbf{Proof}

The proof of this result, which can be executed via an application of the classical Div-Curl Lemma, is given in \cite[Lemma 8.1]{afn}.

$\hfill \Box$

In the limit passage in the Galerkin method in Section \ref{gallim}, we are not able to identify the weak limit of the boundary term $|\phi + u_n|$, cf.\ the convergence \eqref{boundconvn}, in the classical sense. However, we are able to identify it in the sense of Young measures. To this end we make use of the following existence result for Young measures in the setting of $L^1$-spaces on the boundaries of Lipschitz domains.

\begin{satz} \label{youngmeasures}
Let $T>0$ and let $\Omega \subset \mathbb{R}^3$ be a bounded domain of class $C^{0,1}$ with boundary $\Gamma := \partial \Omega$. Let $(u_n)_n$ be a sequence of functions bounded (uniformly with respect to $n$) in $L^1((0,T)\times \Gamma; \mathbb{R}^3)$. Then, there exists a parametrized family $\{\nu_{(t,y)}\}_{(t,y)\in (0,T)\times \Gamma}$ of probability measures $\nu_{(t,y)}$ on $\mathbb{R}^3$, depending measurably on $(t,y)$, such that for any continuous function $\Psi: (0,T) \times \Gamma \times \mathbb{R}^3 \rightarrow \mathbb{R}$ with the property
\begin{align}
\Psi \left(\cdot, \cdot ; u_n(\cdot, \cdot) \right) \rightharpoonup \overline{\Psi (u)}(\cdot, \cdot) \quad \text{in } L^1\left((0,T) \times \Gamma \right) \label{boundconv}
\end{align}

it holds that
\begin{align}
\overline{\Psi (u)}\left(t,y \right) = \int_{\mathbb{R}^3} \Psi \left(t,y;z \right)\ d\nu_{(t,y)}(z) \quad \quad \text{for a.a. } (t,y) \in (0,T)\times \Gamma. \nonumber
\end{align}
\end{satz}

\textbf{Proof}

The idea behind the proof is a transfer of the classical theory of Young measures (see e.g.\ \cite[Theorem 0.10]{fn}) to the setting of $L^1$-spaces on the boundaries of Lipschitz domains. We begin by introducing the notation necessary for the definition of weak convergence in the space $L^1((0,T)\times \Gamma)$. Due to the $C^{0,1}$-regularity of $\Omega$, we find open sets $U_j \subset \mathbb{R}^3$, $j=1,...,J \in \mathbb{N}$, covering a neighbourhood $U \subset \mathbb{R}^3$ of $\Gamma$ and for each $j=1,...,J$ a local coordinate system $(y',y_3)$, an open set $U_j' \subset \mathbb{R}^2$ and a bi-Lipschitz function $h_j \in C^{0,1}(U_j')$ such that
\begin{align}
U_j \bigcap \Gamma = \left\lbrace \left(y',h_j\left(y' \right) \right):\ y' \in U_j' \right\rbrace. \nonumber
\end{align}

We further find a partition of unity $\left\lbrace \Phi_j \right\rbrace_{j=1}^J$ in $U$,
\begin{align}
\Phi_j \in \mathcal{D}\left(U_j \right),\quad \quad 0 \leq \Phi_j \leq 1 \quad \text{in } U_j,\quad \quad \sum_{j=1}^J \Phi_j = 1 \quad \text{in } U. \nonumber
\end{align}

The weak convergence \eqref{boundconv} can then be expressed as
\begin{align}
&\sum_{j=1}^J \int_0^T \int_{U_j'} \Phi_j \left(y',h_j\left(y'\right) \right) \phi \left(t,y',h_j\left(y'\right) \right) \Psi \left(t,y',h_j\left(y'\right); u_n \left(t,y',h_j\left(y'\right) \right) \right) \sqrt{ 1 + \left| Dh_j\left(y'\right) \right|^2}\ dy'dt \nonumber \\
\rightarrow & \sum_{j=1}^J \int_0^T \int_{U_j'} \Phi_j \left(y',h_j\left(y'\right) \right) \phi \left(t,y',h_j\left(y'\right) \right) \overline{\Psi \left( u \right)} \left(t,y',h_j \left(y'\right) \right)  \sqrt{ 1 + \left| Dh_j\left(y'\right) \right|^2}\ dy'dt \label{defofweakconv}
\end{align}

for all $\phi \in L^\infty ((0,T)\times \Gamma)$. We now fix some arbitrary point $x \in \Gamma$. Since the weak convergence in $L^1((0,T)\times \Gamma)$ is in fact independent of the specific choice of the sets $U_j$ and the functions $\Phi_j$ we may assume, without loss of generality, that $\Phi_j = 1$ in a small neighbourhood $U_x \subset U_j$ of $x$ for some $j=1,...,J$. Hence, for any $\phi \in L^\infty ((0,T) \times (U_x \bigcap \Gamma))$, the convergence \eqref{defofweakconv} reduces to
\begin{align}
&\int_0^T \int_{U_j'} \phi \left(t,y',h_j\left(y'\right) \right) \Psi \left(t,y',h_j\left(y'\right); u_n \left(t,y',h_j \left(y'\right) \right) \right) \sqrt{ 1 + \left| Dh_j\left(y'\right) \right|^2}\ dy'dt \nonumber \\
\rightarrow & \int_0^T \int_{U_j'} \phi \left(t,y',h_j\left(y'\right) \right) \overline{\Psi \left( u \right)} \left(t,y',h_j \left(y'\right) \right) \sqrt{ 1 + \left| Dh_j\left(y'\right) \right|^2}\ dy'dt. \nonumber
\end{align}

Since $\sqrt{1 + |Dh_j|^2}$ is bounded away from zero, it then further follows that
\begin{align}
\Psi' \left(\cdot, \cdot; u_n'(\cdot, \cdot) \right) \rightharpoonup \overline{\Psi (u)}(\cdot, \cdot, h_j(\cdot)) \quad \quad \text{in } L^1\left((0,T) \times U_x' \right), \label{convinU'}
\end{align}

where
\begin{align}
\Psi'\left(t,y';z \right) := \Psi \left(t,y',h_j\left(y' \right);z \right),\quad \quad u_n'\left(t,y'\right) := u_n \left(t,y',h_j\left(y'\right) \right) \quad \quad \quad \forall t \in (0,T),\ y' \in U_x',\ z \in \mathbb{R}^3 \nonumber
\end{align}

and
\begin{align}
U_x' = \left\lbrace y' \in U_j':\ \left(y',h_j'\left(y' \right)\right) \in U_x \bigcap \Gamma \right\rbrace. \nonumber
\end{align}

At this point, we are in the position to apply the classical theory of Young measure's (see for example \cite[Theorem 0.10]{fn}) for the identification of the limit function in the convergence \eqref{convinU'}. This guarantees the existence of a parametrized family $\{\nu_{(t,y')}\}_{(t,y')\in (0,T)\times U_x'}$ of probability measures $\nu_{(t,y')}$ on $\mathbb{R}^3$, depending measurably on $(t,y')$ and independent of $\Psi$, such that
\begin{align}
\overline{\Psi (u)}\left(t,y',h_j\left(y' \right) \right) &= \int_{\mathbb{R}^3} \Psi' \left(t,y';z \right)\ d\nu_{(t,y')}(z)\nonumber \\
&= \int_{\mathbb{R}^3} \Psi \left(t,y',h_j\left(y'\right);z \right)\ d\nu_{(t,y')}(z) \quad \quad \text{for a.a. } (t,y') \in (0,T)\times U_x'. \nonumber
\end{align}

Since any $y \in U_x \bigcap \Gamma$ can be identified with some $y' \in U_x'$ via the relation $y = (y',h_j(y'))$, this means that
\begin{align}
\overline{\Psi (u)}\left(t,y \right) = \int_{\mathbb{R}^3} \Psi \left(t,y;z \right)\ d\nu_{(t,y)}(z) \quad \quad \text{for a.a. } (t,y) \in (0,T)\times \left( U_x \bigcap \Gamma \right), \label{desidentity}
\end{align}

where we defined $\nu_{(t,y)} := \nu_{(t,y')}$. Due to the arbitrary choice of $x \in \Gamma$, the identity \eqref{desidentity} actually holds true for almost all $(t,y) \in (0,T)\times \Gamma$ and thus we conclude the proof.

$\hfill \Box$

\subsection* {Acknowledgments:}
  The work of \v S. N.  and J. S. was supported by Praemium Academiae of \v S. Ne\v casov\' a and by the Czech Science Foundation (GA\v CR) through project GA22-01591S. The Institute of Mathematics, Czech Academy of Sciences, is supported by RVO:67985840.

\end{document}